\documentclass{gtart_h}  

\def\ifplaintex{\expandafter\ifx\csname documentclass\endcsname\relax}

\ifplaintex 
\else
\fi

\expandafter\ifx\csname beginpicture\endcsname\relax
\expandafter\ifx\csname documentclass\endcsname\relax
\input pictex \else
\input prepictex \input pictex \input postpictex \fi\fi

\def\gt{{\mathsurround=0pt\it $\cal G\mskip-2mu$eometry \&\ 
$\cal T\!\!$opology}}        %  journal title in recommended style

\def\gtp{{\mathsurround=0pt\it $\cal G\mskip-2mu$eometry \&\ 
$\cal T\!\!$opology $\cal P\!$ublications}}  % GT publications

\def\lognumber#1{\def\thelognumber{#1}}
\def\volumenumber#1{\def\thevolumenumber{#1}}
\def\papernumber#1{\def\thepapernumber{#1}}
\def\volumeyear#1{\def\thevolumeyear{#1}}

\def\pagenumbers#1#2{\def\startpage{#1}\def\finishpage{#2}}
\def\published#1{\def\publishdate{#1}}
\def\proposed#1{\def\theproposer{#1}}
\def\seconded#1{\def\theseconders{#1}}
\def\received#1{\def\receiveddate{#1}}
\def\revised#1{\def\reviseddate{#1}}
\def\accepted#1{\def\accepteddate{#1}}

\long\def\asciiabstract#1{\long\def\theasciiabstract{#1}}

\let\\\par\let\thelognumber\relax
\let\thevolumenumber\relax\let\thepapernumber\relax
\let\thevolumeyear\relax\let\thesamplenumber\relax\let\startpage\relax
\let\finishpage\relax\let\publishdate\relax\let\receiveddate\relax
\let\reviseddate\relax\let\accepteddate\relax\let\theasciititle\relax
\let\theasciiauthors\relax
\let\theasciiabstract\relax
\let\theasciiemail\relax\let\theshortauthors\relax\let\theshorttitle\relax

\long\def\maketitlep{   % start of definition of \maketitlep

\count0=\startpage

\gt\hfill      %   Journal title (top left) 
\beginpicture
\setcoordinatesystem units <0.33truein, 0.33truein> point at 2.2 0.9
\setplotsymbol ({$\cal G$})
\plotsymbolspacing=9truept
\circulararc 315 degrees from 0 1 center at 0 0
\setplotsymbol ({$\cal T$})
\circulararc 315 degrees from 1 -1 center at 1 0
\endpicture
\break
{\small\ifx\thesamplenumber\relax % sample?  
Volume \else Sample
\fi\thevolumenumber\ (\thevolumeyear)
\startpage--\finishpage\nl
Published: \publishdate}
\vglue 0.5truein plus 0.4fil minus 0.1truein

{\parskip=0pt\leftskip 0pt plus 1fil\def\\{\par\smallskip}{\ifplaintex\large
\else\Large\fi\bf\thetitle}\par\medskip}   

\vglue 0pt plus 0.1fil 

{\parskip=0pt\leftskip 0pt plus 1fil\def\\{\par}{\sc\theauthors}
\par\medskip}

\vglue 0pt plus 0.1fil 

{\small\parskip=0pt\let\newline\\
{\leftskip 0pt plus 1fil\def\\{\par}{\sl\theaddress}\par}
\expandafter\ifx\theemail\relax    % email address?
\relax\else\vglue 5pt plus 0.02fil minus 2pt\def\\{\stdspace{\rm 
and}\stdspace} 
\cl{Email:\stdspace\tt\theemail}\fi
\ifx\theurl\relax                  % URL given?
\relax\else\vglue 5pt plus 0.02fil minus 2pt\def\\{\stdspace{\rm 
and}\stdspace}
\cl{URL:\stdspace\tt\theurl}\fi\par}

\vglue 7pt plus 0.3fil minus 3pt

{\bf Abstract}
\vglue 5pt plus 0.1fil minus 2pt

\theabstract

\vglue 7pt plus 0.3fil minus 3pt

{\bf AMS Classification numbers}\quad Primary:\quad \theprimaryclass

Secondary:\quad \thesecondaryclass

\vglue 5pt plus 0.3fil minus 2pt

{\bf Keywords:}\quad \thekeywords

\vglue 10pt plus 0.5fil minus 5pt

{\small  Proposed: \theproposer\hfill Received: \receiveddate\nl
Seconded: \theseconders\hfill 
\ifx\reviseddate\relax                         % paper revised?
Accepted: \accepteddate                        % no
\else
Revised: \reviseddate                          % yes
\fi}
\eject
}       %  end of definition of \maketitlep

\let\maketitlepage\maketitlep
\let\maketitle\maketitlepage

\font\phead=cmsl9 scaled 950
\font=cmsl9 scaled 1050
\font\pnum=cmbx10 scaled 913
\font\lnum=cmbx10 
\font\pfoot=cmsl9 scaled 950
\font=cmsl9 scaled 1050
\ifplaintex
\headline{\vbox to 0pt{\vskip -4.5mm\line{\small\phead\ifnum
\count0=\startpage ISSN 1364-0380 (on line)
1465-3060 (printed) \hfill {\pnum\folio}\else\ifodd\count0\def\\{ }% 
\ifx\theshorttitle\relax\thetitle\else\theshorttitle\fi\hfill{\pnum\folio}
\else\def\\{ and }{\pnum\folio}\hfill\ifx\theshortauthors\relax\theauthors
\else\theshortauthors\fi\fi\fi}\vss}}
\footline{\vbox to 0pt{\vglue 0mm\line{\small\pfoot\ifnum\count0=\startpage
\copyright\ \gtp\hfill\else
\gt, Volume \thevolumenumber\ (\thevolumeyear)\hfill\fi}\vss
}}
\else
\makeatletter
\def\@oddhead{{\small\lhead\ifnum\count0=\startpage ISSN 1364-0380 (on line)
1465-3060 (printed) \hfill {\lnum\number\count0}\else\ifodd\count0
\def\\{ }\ifx\theshorttitle\relax \thetitle \else\theshorttitle\fi\hfill
{\lnum\number\count0}\else\def\\{ and }{\lnum\number\count0}
\hfill\ifx\theshortauthors\relax 
\theauthors\else\theshortauthors\fi\fi\fi}}\def\@evenhead{\@oddhead}
\def\@oddfoot{\small\lfoot\ifnum\count0=\startpage\copyright\ \gtp\hfill\else
\gt, Volume \thevolumenumber\ (\thevolumeyear)\hfill\fi}
\def\@evenfoot{\@oddfoot}
\makeatother
\fi

\newwrite\gtoutfile
\long\gdef\makeheadfile{  %%% start of definition of \makeheadfile
{\def\\{, }\def\s{ }
\immediate\openout\gtoutfile head.xxx
\immediate\write\gtoutfile{Proxy-for: \ifx\theasciiauthors\relax
\theauthors\else\theasciiauthors\fi\s<\ifx\theasciiemail\relax\theemail\else\theasciiemail\fi>}
\immediate\write\gtoutfile{\noexpand\\}
\immediate\write\gtoutfile{Authors: \ifx\theasciiauthors\relax
\theauthors\else\theasciiauthors\fi}
{\def\\{ }\immediate\write\gtoutfile{Title: \ifx\theasciititle\relax
\thetitle\else\theasciititle\fi}}
\immediate\write\gtoutfile{Subj-class: GT or SG or MG etc}
\immediate\write\gtoutfile{MSC-class: \theprimaryclass\ifx\thesecondaryclass\relax\else, \thesecondaryclass\fi}
\immediate\write\gtoutfile{Journal-ref: Geom. Topol. \thevolumenumber
(\thevolumeyear) \startpage-\finishpage}
\immediate\write\gtoutfile{Comments: Published by Geometry and Topology at}
\immediate\write\gtoutfile{\s\s http://www.maths.warwick.ac.uk/gt/GTVol\thevolumenumber/paper\thepapernumber.abs.html}
\immediate\write\gtoutfile{\noexpand\\}
\immediate\write\gtoutfile{}
\ifx\theasciiabstract\relax
\immediate\write\gtoutfile{\theabstract}\else
\immediate\write\gtoutfile{\theasciiabstract}\fi
\immediate\write\gtoutfile{}
\immediate\write\gtoutfile{\noexpand\\}
\immediate\write\gtoutfile{}
\immediate\closeout\gtoutfile}}  %%% end of definition of \makeheadfile

\def\maketitlepage{\maketitlep\makeheadfile}
\let\maketitle\maketitlepage

\lognumber{355}
\received{22 August 2003}
\volumenumber{9}\papernumber{26}\volumeyear{2005}
\pagenumbers{1147}{1185}   
\revised{9 March 2005}
\published{14 June 2005}
\accepted{8 June 2005}
\proposed{Martin Bridson}
\seconded{Dieter Kotschick, Benson Farb}

\usepackage{amsmath,amssymb}
\usepackage{color}
\usepackage{graphicx}

\def\R{\mathbb{R}}
\def\Z{\mathbb{Z}}

\def\N{\mathbb{N}}
\def\Q{\mathbb{Q}}

\newtheorem{theorem}{Theorem}[section]
\newtheorem{lemma}[theorem]{Lemma}
\newtheorem{corollary}[theorem]{Corollary}
\newtheorem{proposition}[theorem]{Proposition}
\newtheorem{question}[theorem]{Question}

\newtheorem{claim}[theorem]{Claim}

\theoremstyle{remark}
\newtheorem{remark}[theorem]{Remark}
\newtheorem{definition}[theorem]{Definition}
\newtheorem{example}[theorem]{Example}

\newcommand{\cobf}{\ensuremath{\|\delta f\|}}

\newcommand{\epf}[1]{\ensuremath{\epsilon_{f,#1}}}
\newcommand{\Zhalf}{\ensuremath{\Z+\frac{1}{2}}}

\begin{document}
\title{Geometry of pseudocharacters}
\author{Jason Fox Manning}
\address{Mathematics 253--37, California Institute of
  Technology\\Pasadena, CA  91125, USA}
\email{manning@caltech.edu}

\begin{abstract}
If $G$ is a group, a pseudocharacter $f\co G\to\R$ is a function which is
``almost'' a homomorphism.  If $G$ admits a nontrivial
pseudocharacter $f$, we define the space of ends of $G$ relative to $f$ and show
that if the space of ends is complicated enough, then $G$ contains a nonabelian
free group.  We also construct a quasi-action by $G$ on a tree whose space of
ends contains the space of ends of $G$ relative to $f$.  
This construction gives rise to
examples of ``exotic'' quasi-actions on trees.  
\end{abstract} 

\asciiabstract{%
If G is a group, a pseudocharacter f: G-->R is a function which is
"almost" a homomorphism.  If G admits a nontrivial pseudocharacter f,
we define the space of ends of G relative to f and show that if the
space of ends is complicated enough, then G contains a nonabelian free
group.  We also construct a quasi-action by G on a tree whose space of
ends contains the space of ends of G relative to f.  This construction
gives rise to examples of "exotic" quasi-actions on trees.}

%geometric group theory
%trees, cohomology of groups


\maketitle

\section{Introduction}
Let $G$ be a finitely presented group.  Following the 
terminology in
\cite{grigorchuk:bounded} a \emph{quasicharacter}\footnote{Several
  competing terminologies exist 
  in the 
  literature.  What we call \emph{quasicharacters}
  are often called
  \emph{quasi-morphisms} \cite{kotschick:what,bavard:longeurstable}
  or \emph{quasi-homomorphisms}
  \cite{bestvinafujiwara:mcg,kotschick:quasi}; what we call
  pseudocharacters are then called
  \emph{homogeneous} quasi-morphisms or quasi-homomorphisms.
  The term ``pseudocharacter'' seems to originate in the papers of
  Fa{\u\i}ziev and Shtern (eg \cite{faiziev:1987,shtern:rademacher}).
  The interest in functions (into normed groups) which are ``almost
  homomorphisms'' goes back at least to Ulam \cite{ulam:problems}.}  
$f$  on a group $G$ is a real valued function on $G$ which 
is a ``coarse homomorphism'' in the sense that the quantity
$f(xy)-f(x)-f(y)$ is 
bounded.  The quasicharacter $f$ is a
\emph{pseudocharacter}
if in addition
$f$ is a homomorphism on each cyclic subgroup of $G$. Any
quasicharacter differs from some pseudocharacter by a bounded function
on $G$.  Brooks \cite{brooks:remarks} gave the first examples where
this pseudocharacter could not be chosen to be a homomorphism.
Brooks' examples were on a free group, but his methods have since been
generalized to give many examples of such ``nontrivial''
quasicharacters on groups
\cite{brooksseries,epsteinfujiwara,fujiwara:gromovhyperbolic,bestvinafujiwara:mcg}.
Interesting applications of such 
existence results can be found in \cite{bavard:longeurstable} and
\cite{kotschick:quasi}.  For additional information on pseudocharacters
we refer the reader to
\cite{kotschick:what,bavard:longeurstable,grigorchuk:bounded} and the
references therein.

Our study of the geometry
of pseudocharacters is partly inspired by the work of Calegari in 
\cite{calegari:cochain}, where it is
shown that if a pseudocharacter on the fundamental group of a $3$--manifold
satisfies some simple geometric hypotheses, then the $3$--manifold must satisfy
a weak form of the Geometrization Conjecture.  The current paper is partially an
attempt to understand what happens when we are given a pseudocharacter which
does \emph{not} satisfy Calegari's hypotheses.  Despite this motivation, we make
no assumptions on the groups considered in this paper, except requiring that
they be finitely presented.

Here is a brief outline of the paper.  In Section \ref{section:definition} we
define the set of ends of a group relative to a pseudocharacter and establish
some basic properties.  In Section \ref{section:action} we make the distinction
between uniform, unipotent and bushy pseudocharacters.  As a kind of warm-up for
the next section we show that a group admitting a bushy pseudocharacter contains
a nonabelian free subgroup.  In Section \ref{section:quasiaction}, we prove the
main theorem:

\medskip
\textbf{Theorem \ref{th:pseudoquasibushy}}\qua
\emph{If $f\co G\to \R$ is a pseudocharacter which is not uniform,
then $G$ admits a cobounded quasi-action on a bushy tree.}
\medskip

We obtain this quasi-action via an isometric action on a space quasi-isometric
to a tree.  A (possibly new) characterization of such spaces is given in
Theorem \ref{th:bottleneck}.  Finally
we use a result of Bestvina and Fujiwara to obtain:

\medskip
\textbf{Theorem \ref{th:infdim}}\qua
\emph{If $G$ admits a single bushy pseudocharacter, then $H_b^2(G;\R)$
and the space of pseudocharacters on $G$ both have dimension equal
to $|\R|$.}
\medskip

Section \ref{section:examples} contains some examples involving negatively curved
$3$--manifolds.  Specifically we show that all but finitely many Dehn surgeries
on the figure eight knot have fundamental groups admitting bushy
pseudocharacters.  This gives the following:

\medskip
\noindent
\textbf{Corollary \ref{cor:counter}}\qua
\emph{There are infinitely many closed $3$--manifold groups which quasi-act
coboundedly on bushy trees but which admit no nontrivial isometric action
on any $\R$--tree.}
\medskip

This corollary can be thought of as an ``irrigidity'' result about quasi-actions
on bushy trees, to be contrasted with the rigidity result of
\cite{msw:quasiactI} about quasi-actions on \emph{bounded valence} bushy trees.

\rk{Acknowledgements}

I would like to thank my PhD advisor Daryl Cooper for his guidance,
support and inspiration.  This work was partially supported by the
NSF, and by a UCSB Graduate Division Dissertation Fellowship.  Thanks
also to the Oxford Mathematical Institute for hospitality while part
of this work was being done, and to the referee for several useful
comments.

\section{Definition of $E(f)$}\label{section:definition}

\subsection{$E(f,S)$ as a set}
\begin{definition}
If $G$ is a finitely presented group, then $f\co G\to \R$ is a
\emph{pseudocharacter} if it has the following properties:
\begin{itemize}
\item $f(\alpha^n)= n f(\alpha)$ for all $\alpha \in G$, $n \in \Z$.
\item $\delta f (\alpha, \beta) = f(\alpha)+f(\beta)- f(\alpha \beta)$ is
bounded independent of $\alpha$ and $\beta$.  We use \cobf\ to denote
the smallest nonnegative $C$ so that $|\delta f(\alpha,\beta)|\leq C$ for all
$\alpha$, $\beta$ in $G$.  
\end{itemize}
\end{definition}

We fix a group $G$ and a pseudocharacter $f\co G\to\R$.  In order to better
understand $f$, we will define a $G$--set $E(f)$, which may be thought of as 
the set of ends of $G$ relative to $f$.

Let $S$ be a finite generating set for $G$.  For simplicity, we assume that
there is a presentation $G=\langle S,R\rangle$ which is triangular, that
is, every
word in $R$ has length three.  It is not hard to show that any finitely
presented group admits a finite triangular presentation.  We will first define
a set $E(f,S)$ and then show it is independent of $S$.

\begin{definition}
If $S$ is a generating set for $G$, let
$\epf{S} = \sup_{s\in S}\{|f(s)|\} + \cobf$.
\end{definition}

If $\Gamma(G,S)$ is the Cayley graph associated to the generating set $S$, we
can extend $f$ affinely over the edges of $\Gamma(G,S)$.  Notice that
\epf{S}\ gives an upper bound on the absolute value of the difference between 
$f(x)$ and $f(y)$ in terms of the distance in $\Gamma(G,S)$ between $x$ and
$y$.  Namely, $|f(x)-f(y)|\leq \epf{S} d_{\Gamma}(x,y)$, if $d_\Gamma$ is the
distance in the Cayley graph.  

\begin{definition}\label{def:sigma}
If $\phi\co\R_+\to\Gamma(G,S)$ is an infinite ray, we define the \emph{sign} of
$\phi$ to be
  \[\sigma_f(\phi) = \left\{\begin{array}{rl}
		+1 & \mathrm{if}\ \lim_{t\to\infty}f\circ\phi(t) = \infty\\
		-1 & \mathrm{if}\ \lim_{t\to\infty}f\circ\phi(t) = -\infty\\
		0 & \mathrm{otherwise.} 
  \end{array}\right.\]
If $f$ is understood, we simply write $\sigma(\phi)$.  If $w$ is some infinite
word in the generators $S$, there is a path $\phi_w\co\R_+\to\Gamma(G,S)$ 
beginning at $1$ and realizing the word.  
We define  $\sigma(w) =\sigma(\phi_w)$.  If $g$ is a group element, we let
$\sigma(g)$ be the sign of $f(g)$.  Notice that if we pick a word $w$
representing $g$ then $\sigma(www\ldots)=\sigma(w^\infty)=\sigma(g)$.
\end{definition}

We give two equivalent definitions of $E(f,S)$. 
\begin{definition}[\rm(Version 1)]\label{def:v1}
 \[E(f,S) = \bigl\{\phi\co \R_+\to \Gamma(G,S)\ \text{continuous}\ |\
 \sigma(\phi)\in\{+1,-1\}\bigr\}/\sim\]
We will say $\phi_1 \sim_C \phi_2$ if
  $\sigma(\phi_1)=\sigma(\phi_2)$ and
for all $D$ with $\sigma(\phi_1)D>C$ 
there is a path $\delta\co[0,1]\to \Gamma(G,S)$ such that:
\begin{itemize}
\item $\delta(0) \in \phi_1(\R_+)$.
\item $\delta(1) \in \phi_2(\R_+)$.
\item $|f \circ \delta (t) - D| \leq C$ for all $t \in [0,1]$.
\end{itemize}
The path $\delta$ will be referred to as a \emph{connecting path}.
We say $\phi_1\sim \phi_2$ if $\phi_1 \sim_C \phi_2$ for some $C$.  This is an
equivalence relation.
\end{definition}
\begin{definition}[\rm(Version 2)]\label{def:v2}
\[E(f,S)=\bigl\{w=w_1 w_2 w_3 \ldots\ |\ w_i \in S\cup S^{-1}
\ \forall i\in\mathbb{Z}\ \text{and}\ \sigma(w)\in\{+1,-1\}\bigr\}/\sim\]
We say $w=w_1 w_2\ldots \sim_C v=v_1 v_2 \ldots$ if
$\sigma(w)=\sigma(v)$
and
 for all $D$ with 
$\sigma(w)D>C$ there is a
word $d = d_1\ldots d_n$ in the letters $S\cup S^{-1}$ such that:
\begin{itemize}
\item $w_p d = v_p$ in $G$ 
for some prefix $w_p$ of $w$ and some prefix $v_p$ of $v$.
\item $|f(w_p d_p) - D| \leq C$ for all prefixes $d_p$ of $d$.
\end{itemize}
The word $d$ will be referred to as a \emph{connecting word}.
We say $w\sim v$ if $w\sim_C v$ for some $C$.  Again, this is an equivalence
relation.
\end{definition}

\begin{lemma}
There is a canonical bijection between $E(f,S)$ (version 1) and $E(f,S)$
(version 2).
\end{lemma}
\begin{proof}
  Let $E_1$
be the set described in Definition \ref{def:v1}, and let $E_2$ be the set
described in Definition \ref{def:v2}.  If $[\phi]\in E_1$, note that
altering $\phi$ on any compact subset of $\R_+$ does not change its
equivalence class.  Thus we may assume $\phi(0)$ is the identity
element of $G$.  The equivalence class of $\phi$ is also left
undisturbed by proper homotopies and arbitrary reparameterizations,
and so we may assume that $\phi$ is a unit speed path with no
backtracking.  Such a representative $\phi$ determines an
infinite word in the generators $w=w_1w_2w_3\ldots$ with each
$w_i=\phi(i-1)^{-1}\phi(i)\in S\cup S^{-1}$.  The reader may check
that this recipe for building a word from a path gives a well-defined
map from $E_1$ to $E_2$, and that this map is a bijection.
\end{proof}

Roughly speaking, if one thinks of the Cayley graph as divided up into thick
slabs of group elements all sent to roughly the same real values, and two paths
pass to infinity through the same slabs, we identify the paths.

\begin{definition} If $[\phi] \in E(f,S)$ and  $\sigma(\phi)=1$
 we say that $[\phi]$ is \emph{positive}.  If $\sigma(\phi)=-1$
 we say that $[\phi]$ is \emph{negative}.  We denote by $E(f,S)^+$ the set of
positive elements and by $E(f,S)^-$ the set of negative elements. 
\end{definition}

\begin{remark}\label{remark:scale}
Notice that neither Definition \ref{def:v1} nor Definition
\ref{def:v2} depend on anything but $S$ and $f$ up to multiplication by a
nonzero real number.  Consequently, $E(f,S)^+$ and $E(f,S)^-$ only depend on 
$S$ and $f$ up to
multiplication by a positive real number.  Thus in what follows we do not
hesitate to scale $f$ by a positive real number whenever convenient.
\end{remark}

\begin{example}  If $G$ is a free abelian group generated by $S$
and 
$f$ is any nontrivial
homomorphism, then $E(f,S)$ contains precisely two elements, 
one positive and one
negative.  Thus $|E(f,S)^+|=|E(f,S)^-|=1$.
\end{example}

\begin{lemma}\label{lemma:inducedaction}
The action of $G$ on the Cayley graph induces an action on $E(f,S)$.  
\end{lemma}
\begin{proof}
This is clearest looking at Definition \ref{def:v1}.  One must only check that
$g\phi\sim g\phi'$ if $\phi \sim \phi'$, where $\phi$ and $\phi'$ are
infinite rays in the Cayley graph with $\sigma(\phi)$ and $\sigma(\phi')$
nonzero.  Suppose that $\phi\sim\phi'$.  Then $\phi\sim_C\phi'$ for some $C>0$.
Let $C'=C+|f(g)|+\cobf$, and suppose that $\sigma(\phi)D>C'$.  Then in
particular $\sigma(\phi)D>C$, and so there is a connecting path $\delta$ with
$\delta(0)\in \phi(\R_+)$, $\delta(1)\in\phi'(\R_+)$, and satisfying
$|f\circ\delta(t)-D|\leq C$ for all $t\in[0,1]$.  Let $\delta'$ be the same
path, translated by $g$.  Then $\delta'(0)\in g\phi(\R_+)$ and $\delta'(1)\in
g\phi'(\R_+)$.  Since $|f\circ\delta'(t)-(f\circ\delta(t)+f(g))|\leq\cobf$, we
have
$|f\circ\delta'(t)-D|\leq|f\circ\delta(t)+f(g)-D|+\cobf
\leq|f\circ\delta(t)-D|+|f(g)|+\cobf
\leq C+|f(g)|+\cobf = C'$.  Thus if $\phi\sim_C\phi'$, then 
$g\phi\sim_{C'} g\phi'$, and so $g\phi\sim g\phi'$.
\end{proof}

\subsection{Topology of $E(f,S)$}\label{section:topology}
We next describe the topology on $E(f,S)$, by describing a basis of open sets in
terms of Definition \ref{def:v1} above.  
\begin{definition}\label{def:topology}
Let
$I$ be some interval in $\R$ of diameter bigger than \epf{S}.  Let $B$ be some
component of $f^{-1}(I) \subset \Gamma(G,S)$.  Let $C$ be some connected
component of $\Gamma(G,S) \setminus B$.  We define
\[ U_{B,C} = \{[\phi]\in E(f,S)\ |\ \mathrm{image}(\phi)\subset C\}\]
\end{definition}
We make $E(f,S)$ a topological space with the collection of all such 
$U_{B,C}$ as a basis.

Note that if $C$ and $C'$ are distinct components of the complement of
$B$, then $U_{B,C}\cap U_{B,C'}$ is empty.

It turns out that $E(f,S)$ is Hausdorff and totally disconnected.  In fact, we
will show the following:
\begin{proposition}\label{prop:tree}
There is a simplicial tree $T$ and a map 
$i\co E(f,S) \to \partial T$ which is 
a homeomorphism onto its image. (By $\partial T$ we mean the Gromov boundary of
$T$.)
\end{proposition}

As noted in Remark \ref{remark:scale}, $E(f,S)$ is unchanged if we scale $f$ by
a nonzero real number.  By scaling appropriately, we can ensure that 
$\epsilon_{f,S} < \frac{1}{4}$, 
and that $f^{-1}\bigl(\Z+\frac{1}{2}\bigr)$ contains no element of $G$.
These assumptions will be made for the rest of the section.

\begin{figure}[ht!]
\begin{center}
\input{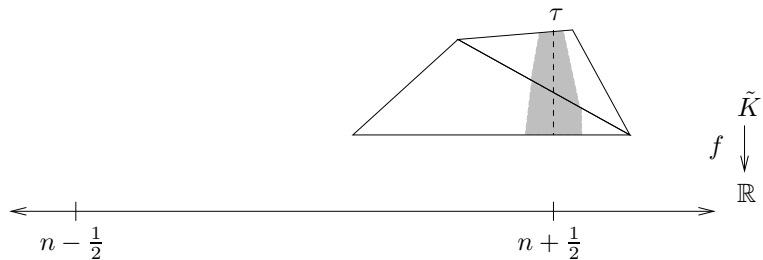}
\caption{Each track $\tau$ has a product neighborhood.}
\label{figure:track}
\end{center}
\end{figure}

Since we chose a triangular presentation of $G$, we may equivariantly add
$2$--simplices to $\Gamma(G,S)$ to obtain a simply connected $2$--complex
$\widetilde{K}$
corresponding to the presentation we started with.  $G$ acts on
$\widetilde{K}$ with
quotient $K$, where $K$ is a one-vertex $2$--complex with one edge for each
element of $S$, and one triangular cell attached
for each relation in our presentation.
The function $f$ may be
extended affinely over the $2$--simplices of $\widetilde{K}$ to give a function
$f\co\widetilde{K}\to\R$.  
Since $\tau = f^{-1}\bigl(\Z+\frac{1}{2}\bigr)$ misses the 
$0$--skeleton of $\widetilde{K}$, and since we have scaled $f$ so that
$|f(v)-f(w)|\leq \epsilon_{f,S}<\frac{1}{4}$ whenever $v$ and $w$ are endpoints of the same edge,
$f^{-1}\bigl(\Z+\frac{1}{2}\bigr)$ intersects each $2$--cell either not at all or in a
single normal arc.  Thus
$\tau$ is a union of possibly infinite tracks in
$\widetilde{K}$.
Each such track $\tau$ separates $\widetilde{K}$ into two components,
and has a product neighborhood 
$\eta(\tau) = \tau \times \bigl(-\frac{1}{2},\frac{1}{2}\bigr)$ in the
complement of the $0$--skeleton of $\widetilde{K}$ (see Figure \ref{figure:track}).
As there may be vertices arbitrarily close to $\tau$, the topological product
neighborhood $\eta(\tau)$ must be allowed to vary in width from $2$--cell to
$2$--cell and is not necessarily a component of $f^{-1}(I)$ for any interval $I$.

We obtain a quotient space $T$
of $\widetilde{K}$ by smashing each component of $\eta(\tau)$ to an interval and 
each component of the complement of $\eta(\tau)$ to a point. Let
$\pi\co\widetilde{K}\to T$ be the quotient map.
 Clearly $T$ is a
simplicial graph.  Since the preimage of each point under $\pi$
is connected and $\widetilde{K}$ is simply connected, $T$ must be simply
connected.  In particular, $T$ is a tree.  We refer to the preimage of an open
edge of $T$ as an \emph{edge space} of $\widetilde{K}$ and to the preimage of a
vertex as a \emph{vertex space} of $\widetilde{K}$.  We define a map 
$\overline{f}\co T\to\R$ as follows.  If $v$ is a vertex of $T$, then the associated
vertex space lies in  $f^{-1}\bigl(n-\frac{1}{2}, n+\frac{1}{2}\bigr)$ for some
$n$; we let
$\overline{f}(v)=n$ and define $\overline{f}$ on edges by extending
affinely.  Note that $|f(x)-\overline{f}\circ\pi (x)|$ is bounded
independently of $x\in\widetilde{K}$.  

\begin{remark}\label{remark:oneskeleton}
Because $f^{-1}\bigl(\Z+\frac{1}{2}\bigr)$ intersects any edge of 
$\widetilde{K}^1=\Gamma(G,S)$ in at most a single point
and any $2$--cell in at most a single normal arc, 
each vertex space must contain vertices of $\widetilde{K}$, and any
two vertices contained in the same vertex space are actually connected by a path
in the intersection of $\widetilde{K}^1$ and the vertex space.
  Since therefore 
the components of
$\Gamma(G,S) \setminus f^{-1}\bigl(\Z+\frac{1}{2}\bigr)$ are in one to one correspondence
with the vertex spaces, the particular pattern of $2$--simplices added to
form $\widetilde{K}$ is unimportant to the structure of $T$.  
\end{remark}

We now define a map from $E(f,S)$ to $\partial T$.  Let $v_0\in T$ be the vertex
whose associated vertex space contains the identity element of $G$.
Since $T$ is a
tree, $\partial T$ can be identified with the set of geodesic rays in $T$ 
starting at $v_0$.  Given some element $[\phi] \in E(f,S)$, we will associate
such a geodesic ray.  First notice that we may assume that the image of $\phi$ 
contains
$1\in G$.  Now consider the image of $\pi\circ\phi$ in $T$.  By our
choice of $\phi$, this image contains $v_0$.
\begin{lemma}
Suppose that $[\phi]\in E(f,S)$ and that the image of $\phi$ contains $1$.
Then the  image of $\pi\circ\phi$ contains a unique geodesic ray starting at 
$v_0$.  
\end{lemma}
\begin{proof}
Since
$\lim_{t\to\infty}\overline{f}\circ\pi\circ\phi(t)=\lim_{t\to\infty}f\circ\phi(t)=\pm\infty$,
and 
 $|\overline{f}(v)-\overline{f}(w)|\leq d(v,w)$ for $v$, $w$ in
$T$, $\pi\circ\phi$ must eventually leave any finite diameter part of $T$.  
Let $B_R=\{x\in T\ |\ d(x,\pi\circ\phi(0))< R\}$ for $R>0$.  For any $R\geq 0$
there is some $t$ so that $\pi\circ\phi([t,\infty))$ lies in a single component
$C$ of $T\setminus B_R$.  Let $x_R$ be the point in $C$ closest to
$\pi\circ\phi(0)$.  Then $\gamma\co \R_+\to T$ given by $\gamma(t)=x_t$ is the
unique geodesic ray starting at $\pi\circ\phi(0)=v_0$ which is contained
in the image of $\pi\circ\phi$.
\end{proof}

We define a map $i\co E(f,S)\to \partial T$ by sending $[\phi]$ to this ray.

\begin{lemma}
The map $i$ is well-defined.
\end{lemma}
\begin{proof}
Suppose that $[\phi]=[\phi']$, but images of the paths $\pi\circ\phi$ and 
$\pi\circ\phi'$ contain distinct infinite rays $r$ and $r'$. For clarity, we
assume that $[\phi]$ is positive.  The proof for $[\phi]$ negative is much the
same. 

We may modify $\phi$ and $\phi'$ so that the image of $\pi\circ\phi$ is $r$ and
the image of $\pi\circ\phi'$ is $r'$.  Furthermore, we may adjust $\phi$ and
$\phi'$ so that $r$ intersects $r'$ in a single point, $v$.   
Let $N_1=\overline{f}(v)$. 
Since $\phi\sim\phi'$ there is some $C>0$ so that $\phi\sim_C\phi'$.  
Let $N_2>C$.

Let $N=|N_1|+N_2+1$.  Since $N>C$, there is some $t\geq 0$, some $t'\geq 0$, and
some path
$\delta\co[0,1]\to\Gamma(G,S)$ so that $\delta(0)=\phi(t)$,
$\delta(1)=\phi'(t')$, and $|f\circ\delta(x)-N|\leq C$ for all $x\in[0,1]$.
But this path $\delta$ must necessarily pass through $\pi^{-1}(v)$, and 
so there is some $x\in[0,1]$ so that $f\circ\delta(x)<N_1+\frac{1}{2}$.  
But this implies
that $|f\circ\delta(x)-N|>N_2>C$, a contradiction.
\end{proof}
\begin{lemma}
The map $i$ is injective.
\end{lemma}
\begin{proof}
Suppose $i([\phi])=i([\phi'])$.  Let $\gamma$ be a unit speed geodesic ray in $T$ contained
in $\pi\circ\phi(\R_+)\cap \pi\circ\phi'(\R_+)$.  We have
$\lim_{t\to\infty}\overline{f}\gamma(t)=\lim_{t\to\infty}f\circ\phi(t)=\lim_{t\to\infty}f\circ\phi'(t)$,
so $\sigma(\phi)=\sigma(\phi')$.  As in the last lemma, we assume for clarity
that
$\phi$ and $\phi'$ are positive.
 
We may assume that
$\gamma(0)$ is a vertex of $T$, and so $\gamma(k)$ is a vertex of $T$
for $k$ any nonnegative integer.  Let $N = \overline{f}(\gamma(0))$.  By
truncating $\gamma$ if necessary, we may assume that $N$ is positive and is the
smallest value taken by $\overline{f}\circ\gamma$.  
For every integer $n$ between $N$ and
$\infty$, there is a vertex $v_n$ in the image of 
$\gamma$ with
$\overline{f}(v_n) = n$.  Both $\phi$ and $\phi'$ must pass through the vertex space
$V_n\subset \widetilde{K}$ associated to $v_n$.  

For each $n\in\Z$, $n\geq N$, 
pick points $x_n$ and
$x_n'$ on the intersections of the paths $\phi$ and $\phi'$ with $V_n$.  Note
that by Remark \ref{remark:oneskeleton}, $x_n$ and $x_n'$ are connected by a
path $\delta_n$ in the intersection of $\widetilde{K}^1$ with $V_n$.  

We claim that $\phi\sim_C\phi'$ for $C= N+2$.  
Suppose $D>C$,
and let $[D]$ be the integer part of $D$.  Since $[D]> N$, there are points
$x_{[D]}$ and $x_{[D]}'$ on the images of $\phi$ and $\phi'$ respectively.
These points are connected by a path $\delta_{[D]}$ which lies entirely in the
intersection of $\widetilde{K}^1$ with $V_{[D]}$.  Since the image of 
$\delta_{[D]}$ lies
entirely inside $V_{[D]}$,  $|f \circ \delta (t) - D| < \frac{3}{2}<C$  for all
$t$.  Thus $\phi\sim_C\phi'$.  In particular, $\phi\sim\phi'$.
\end{proof}

The following two lemmas complete the proof of Proposition \ref{prop:tree}.
\begin{lemma} 
The map $i$ is continuous.
\end{lemma}
\begin{proof}
Recall that the topology on $\partial T$ can be described by a basis
of open sets, as follows:  Let $e$ be any open edge of $T$, and let $T'$
be one of the two components of $T\setminus e$.  There is a natural inclusion of
$\partial T'$ into $\partial T$; the image of $\partial T'$
is a basic open set.  The
topology on $\partial T$ is generated by such sets.  

Let $U = \partial T'$ be an element of the basis described above.  
Let $v$ be the unique vertex in
$(T\setminus T')\cap \overline{e}$, where $\overline{e}$ is the closed edge
whose interior is $e$. 
If $N(v)=\bigl\{x\in T\ |\ d(x,v)<
\frac{1}{2}\bigr\}$, then $B=\pi^{-1}(N(v))$ is a path component of the
preimage
$f^{-1}\bigl(\overline{f}(v)-\frac{1}{2},\overline{f}(v)+\frac{1}{2}\bigr)$, as it
contains the entire vertex space corresponding to $v$ and half of each
edge space adjacent to $v$.
Furthermore, $\pi^{-1}(T')$ lies
entirely in the complement of $B$.  Let $C$ be the component of the complement
of $B$ containing $\pi^{-1}(T')$, and let $U_{B,C}$ be the basic open set in
$E(f,S)$ defined by $C$ (as in Definition \ref{def:topology}). Then clearly
$i^{-1}(U) = U_{B,C}$.  
Thus $i$ is continuous.
\end{proof}

\begin{lemma}
The map $i$ is an open map onto its image, topologized as a subset of $\partial
T$.  
\end{lemma}
\begin{proof}
If $U_{B,C}$ be given as in Definition
\ref{def:topology}, then $\pi(B)$ is some connected subset of $T$.
Thus $\pi(B)$ is contained in some
minimal (possibly infinite) subtree $Z$.  Because $f$ and
$\overline{f}\circ\pi$ are boundedly different from one another, and every
point in $Z$ is distance at most $1$ from $\pi(B)$, any geodesic in $T$
representing an element of $i(E(f))$ must eventually leave $Z$ forever.
If $e$ is any edge with precisely one endpoint in
$Z$, then let $T_e$ be the component of $T \setminus \mathrm{Int}(e)$ which does
not contain $Z$.  

\begin{claim}
$\partial T_e\cap i(E(f,S)) \subset i(U_{B,C})$ or 
$\partial T_e \cap i(U_{B,C})=\emptyset$.
\end{claim}
\begin{proof}
To prove the claim, suppose there is some point $x$ in 
$\partial T_e\cap i(U_{B,C})$ and suppose that 
$y\in \partial T_e\cap i(E(f,S))$.
Both $i^{-1}(x)$ and $i^{-1}(y)$ are represented by paths $\phi_x$ and $\phi_y$
whose images lie entirely in $\pi^{-1}(T_e)$.  But since
$B\subset\pi^{-1}(Z)$, the space $\pi^{-1}(T_e)$ is contained entirely
in a single complementary component of $B$.  Since $x\in i(U_{B,C})$,
this complementary component is $C$, and so $[\phi_y]\in U_{B,C}$.
Therefore $y\in i(U_{B,C})$.
\end{proof}

By the claim, $i(U_{B,C})$ can be expressed as a union of basic open sets in 
$i(E(f,S))$. The lemma follows.
\end{proof}

\begin{corollary}
$E(f,S)$ is metrizable.  In particular, $E(f,S)$ is Hausdorff.
\end{corollary}

\subsection{Invariance under change of generators}

We have been using $f$ to refer both to the pseudocharacter and to its
extension to $\Gamma(G,S)$.  For this subsection we need to deal with distinct
generating sets, so we will temporarily refer to the extension of $f$ to a
particular Cayley graph $\Gamma(G,S)$ as $f_S$.

Let $T$ be another finite triangular generating set for $G$.  We choose
equivariant maps $\tau\co\Gamma(G,S)\to\Gamma(G,T)$ and
$\upsilon\co\Gamma(G,T)\to\Gamma(G,S)$ which are the 
identity on $G$ and send each edge to a constant speed path.  
Let $N$ be the maximum length of the image of a single edge
under $\tau$ or $\upsilon$.  It is not hard to establish that both $\tau$
and $\upsilon$ are continuous $(N,2N+2)$ quasi-isometries which are
quasi-inverses of one another.  In fact, we have
\[
d(x,y)/N-(N+1/N)\leq d(\tau x,\tau y)\leq N d(x,y)
\]
and similar inequalities for $\upsilon$.

We define a map $\overline{\tau}\co E(f,S)\to E(f,T)$ by
$\overline{\tau}([\phi])=[\tau\circ\phi]$, and define $\overline{\upsilon}$ similarly.

\begin{lemma}
The maps $\overline{\tau}$ and $\overline{\upsilon}$ are well-defined.
\end{lemma}
\begin{proof}
Suppose $[\phi] = [\phi']\in E(f,S)$.  
Then $\phi\sim_C\phi'$ for some $C>0$.  If $\delta$ is a connecting path between
$\phi$ and $\phi'$, then $\tau\circ\delta$ gives a connecting path between
$\tau\circ\phi$ and $\tau\circ\phi'$.  The only trouble is that $f_T$ may vary
on $\tau\circ\delta$ more than $f_S$ varies on $\delta$.  On the other hand,
$\tau\circ\delta$ never gets further than $\frac{N}{2}$ (in $\Gamma(G,T)$ 
from the group elements contained in the image of $\delta$, and so if
$\phi\sim_C\phi'$, then $\tau\circ\phi\sim_{C'}\tau\circ\phi'$ for 
$C'=C+\frac{N}{2}\epsilon_{f,T}$. 
  Thus $[\tau\circ\phi]=[\tau\circ\phi']$.  The
proof for $\overline{\upsilon}$ is identical.
\end{proof}

\begin{lemma}
The maps $\overline{\tau}$ and $\overline{\upsilon}$ are  bijections.
\end{lemma}
\begin{proof}
$\overline{\tau}$ and $\overline{\upsilon}$ are inverses of one another.
\end{proof}

\begin{lemma}
The maps $\overline{\tau}$ and $\overline{\upsilon}$ are open.
\end{lemma}
\begin{proof}

Let $U_{B,C}\subset E(f,S)$ be a basic open set. There is some interval
$[a,b]\subset \R$ so that $B$ is a connected component of
$f_S^{-1}[a,b]\subset\Gamma(G,S)$ and $C$ is a component of
$\Gamma(G,S)\setminus B$.  We wish to show that $\overline{\tau}(U_{B,C})$ is open in
$E(f,T)$.  

Since $\tau$ is continuous, $\tau(B)$ is connected.  As no edge of $B$ has image
of length more than $N$, $f_T\circ\tau(B)\subset(a-N\epf{T},b+N\epf{T})$.
Thus $\tau(B)\subset B'$ a connected component of
$f_T^{-1}[a-(N+N^2+2)\epf{T},b+(N+N^2+2)\epf{T}]$.  
We have chosen the constants here so that the 
distance between any point in the complement of $B'$ and any point in $\tau(B)$
is at least $N^2+2$.

We claim that $\overline{\tau}(U_{B,C})$
is a union of open sets of the form $U_{B',C'}$ where $C'$ is a component of the
complement of $B'$.  The claim follows if each such $U_{B',C'}$ is either
contained in or disjoint from $\overline{\tau}(U_{B,C})$. 

Suppose $[\phi_1]$ and $[\phi_2]$ are in $U_{B',C'}$.  Since $\overline{\tau}$ is
onto and removing an initial segment does not change the equivalence class of a
path, we may suppose $\phi_i=\tau\circ\psi_i$, where each $\psi_i$ has image
entirely in the complement of $B$ and each $\phi_i$ has image entirely in $C'$.
Thus there is a path $\delta$ in $C'$ connecting $\phi_1(0)$ to $\phi_2(0)$.
The path $\sigma\circ\delta$ therefore runs from $\psi_1(0)$ to $\psi_2(0)$.

If $\psi_1(0)$ and $\psi_2(0)$ were in different components of
$\Gamma(G,S)\setminus B$, then $\sigma\circ\delta$ would pass through $B$.  But
then $\sigma\circ\delta$ must pass through a vertex $v$ of $B$ 
(since $B$ is not contained in an edge).  Because
\[d(x,y)/N-(N+1/N)\leq d(\sigma x,\sigma y)\]
the distance between the path
$\delta$ and $\tau(v)$ is less than $N^2+1$.  But this contradicts the
assertion that $\delta$ lies entirely outside of $B'$. 

Since $\psi_1$ and $\psi_2$ are infinite paths in the same component of the
complement of $B$, either both $[\psi_1]$ and $[\psi_2]$ are in $U_{B,C}$, or
neither is.  Likewise, either both $[\phi_1]$ and $[\phi_2]$ are in
$\overline{\tau}(U_{B,C})$ or neither is, establishing the claim.

Again, the proof for $\overline{\upsilon}$ is identical.
\end{proof}

\begin{corollary}\label{cor:invt}
$\overline{\tau}$ is a homeomorphism.
\end{corollary}

\section{The action of $G$ on $E(f)$}\label{section:action}
\subsection{Dynamics}
Let $S$ be a generating set for $G$.  The group $G$ acts on $E(f)=E(f,S)$ 
via the action on the Cayley graph (Lemma \ref{lemma:inducedaction}).
\begin{lemma}
$G$ acts on $E(f)$ by homeomorphisms.
\end{lemma}
\begin{proof}
Let $U_{B,C}$ be a basic open set, so that $B$ is a connected component of
$f^{-1}[a,b]\subset \Gamma(G,S)$ and $C$ is a connected component of the
complement of $B$.  Let $g\in G$.
Note that $g(B)\subset f^{-1}[a+f(g)-\cobf,b+f(g)+\cobf]$.
Let $B'$ be the connected component of $f^{-1}[a+f(g)-\cobf,b+f(g)+\cobf]$
containing $g(B)$.  
If $C'$ is a complementary component of $B'$ we wish to claim
that either $U_{B',C'}\subset g(U_{B,C})$ or $U_{B',C'}$ is disjoint from
$g(U_{B,C})$.  It will follow that $g$ acts on $E(f)$ by an open map.  Since
$g^{-1}$ must do likewise, it follows that $G$ acts by homeomorphisms.

To establish the claim, suppose that $[\phi_1]$ and $[\phi_2]$ are in
$U_{B',C'}$.  We may assume that the images of $\phi_1$ and $\phi_2$ lie
entirely in $C'$.  It follows that $g^{-1}\phi_1$ and $g^{-1}\phi_2$ have image
entirely in the same complementary component of $B$, and the claim is
established.
\end{proof}

Now that we have established that $G$ acts by homeomorphisms of $E(f)$ we look
more closely at the dynamics of this action.  
\begin{definition}
If $a$ and $r$ are fixed points of a group element $g$, we say that $a$ is
\emph{attracting} and $r$ is \emph{repelling} if for any neighborhood $U$ of $a$
and any neighborhood $V$ of $r$, we have $g^n(E(f)\setminus V)\subset U$
for all $n$ sufficiently large. 
\end{definition}

\begin{lemma}\label{lemma:dynamics}
If $g\in G$ and $f(g)\neq 0$, then $g$ has exactly two fixed points in $E(f)$,
one attracting and one repelling.
\end{lemma}
\begin{proof}
It is convenient to use Definition \ref{def:v2} here.
Let $w$ be a word in the letters $S\cup S^{-1}$ representing $g$.  Since
$f(g)\neq 0$, the words $[w^\infty]=[www\ldots]$ and
$[\overline{w}^\infty]=[\overline{w}\overline{w}\overline{w}\ldots]$ are elements of $E(f)$.  Both
are clearly fixed by $g$.

Let $U$ be an open set containing $[w^\infty]$, and let $V$ be an open set
containing $[\overline{w}^\infty]$.  Without loss of generality, both $U$ and $V$ are
basic open sets  $U=U_{B,C}$ and $V=V_{D,E}$.

Let
$\gamma$ be the bi-infinite line made by taking the path from $1$ to $g$
described by $w$ and translating it by powers of $g$.  The group 
$\langle g\rangle\subset G$ acts on $\gamma$ as $\Z$ acts on $\R$.  
Figure \ref{figure:attract} shows approximately how all this might look in $G$.

\begin{figure}[ht!]
\begin{center}
\input{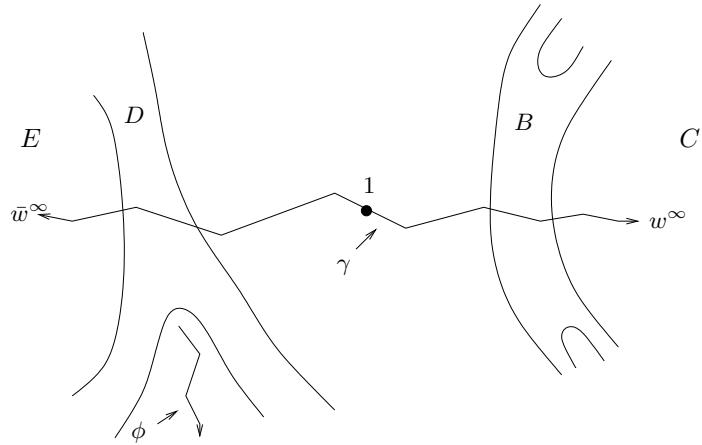}
\caption{Possible arrangement of the barrier spaces}
\label{figure:attract}
\end{center}
\end{figure}

Since $f$ restricted to $\gamma$ is a continuous quasi-isometry,
$D\cap \gamma$ is a compact set.  Thus there is some $N_1$ so that
$g^n(D\cap\gamma)$ is in $C$ for all $n>N_1$.  The barrier spaces $B$ and $C$
are components of the preimages of closed intervals under $f$.  So in particular
there are intervals $[a,b]$ and $[c,d]$ so that $f(B)=[a,b]$ and $f(D)=[c,d]$.
It is easy to check that $f(g^n D)\subset [c+nf(g)-\cobf, d+nf(g)+\cobf]$.  Thus
there is some $N_2$ so that $g^n D\cap B$ is empty for any $n>N_2$.  Let
$N=\max\{N_1,N_2\}$.  It is clear that $g^n(D)\subset C$ for any $n>N$.  

Let $e\in E(f)\setminus V$, and suppose that $n>N$.
We claim that $g^n e\in U_{B,C}$.  
For let $\phi\co G\to\Gamma(G,S)$ be a path so that $[\phi]=e$.  
We may assume that $\phi$ maps
entirely into the complement of $D$, and thus that $g^n\phi$ maps entirely into
the complement of $g^n D$.  But if $g^n[\phi]\notin U_{B,C}$, then eventually
the image of $g^n[\phi]$ leaves the half-space $C$.  Since $g^n D$ is contained
entirely in $C$ this means that $g^n\phi$ maps entirely into the component of
the complement of $g^n D$ which also contains $B$.  In other words, $g^n\phi$
maps into the same component of the complement of $g^n D$ as $\overline{w}^\infty$
does, namely $g^n E$.  But this implies that $e\in U_{D,E}$, a contradiction to
our original choice of $D$ and $E$.
\end{proof}

\subsection{The bushy case}

\begin{definition}
Let $E(f)^+\subset E(f)$ be the set of positive elements of $E(f)$, and let
$E(f)^-$ be the set of negative elements.
\end{definition}
\begin{remark}\label{remark:ginfty}
So long as there exists some $g$ with $f(g)\neq0$, then
$E(f)^+$ and $E(f)^-$ are nonempty.  Indeed, if $w$ is any word
representing $g$, then the infinite word $w^\infty=www\ldots$
determines an element of $E(f)^\pm$, depending on whether $f(g)$ is
positive or negative.  Similarly the infinite word
$\overline{w}^\infty=\overline{w}\overline{w}\overline{w}\ldots$ determines an element of
$E(f)^\mp$.  Neither of these elements actually depends on $w$; we
abuse notation slightly by writing them as  $[g^\infty]$ and
$[g^{-\infty}]$ respectively.   
\end{remark}
\begin{definition}
Let $f$ be a pseudocharacter.  If $|E(f)|=2$ we say $f$ is \emph{uniform}.  If
$|E(f)^+|=1$ or $|E(f)^-|=1$ but $f$ is not uniform, we say $f$ is
\emph{unipotent}.  Otherwise we say that $f$ is \emph{bushy}.
\end{definition}

In \cite{calegari:cochain}, Calegari shows
that if $f$ is uniform and $G$ is the fundamental group of a closed irreducible
$3$--manifold, then $G$
satisfies the Weak Geometrization Conjecture.  Thus one of our goals is to give
information about what happens if $f$ is not uniform.  

\begin{remark}
The terminology here is slightly different from \cite{calegari:cochain}.  What
is here called \emph{uniform} is called \emph{weakly uniform} in
\cite{calegari:cochain}.  For $f$ to be uniform in the sense of
\cite{calegari:cochain}, its coarse level sets must be coarsely simply
connected.
\end{remark}

Following \cite{calegari:cochain}, we
define an \emph{unambiguously positive} element of $G$ to be an
element $g$ with $f(g)>\cobf$.  Note that if $g$ is unambiguously
positive and $h$ is any element of $G$, then
$f(hg)>f(h)$.  If $S$ is any triangular set of generators, we may
alter $S$ so that it contains an unambiguously positive element and is still
triangular.  By Corollary \ref{cor:invt} this has no effect on the $G$--set 
$E(f)$.  It is convenient in what follows to assume that $S$ contains an
unambiguously positive element.  In this case we say that the
generating set $S$ is \emph{unambiguous}.

\begin{lemma}\label{lemma:groupelements}
If $E(f)$ is bushy, then there are group elements $g_1$, $g_2$,  and
$g_3$ so that $[g_1^\infty]\neq[g_2^\infty]$ and
$[g_1^{-\infty}]\neq[g_3^{-\infty}]$ and $f(g_i)>0$ for $i\in\{1,2,3\}$.
\end{lemma}
\begin{proof}
Let $g_1$ be an unambiguously positive element of $S$, where $S$ is some
fixed unambiguous triangular generating set.  Then
$[g_1^\infty]\in E(f)^+$ (see Remark \ref{remark:ginfty} above).  
By assumption there is some
$\phi\co\R_+\to\Gamma(G,S)$ with $[\phi]\in E(f)^+$ but
$[\phi]\neq[g_1^\infty]$.  We may assume that $\phi(0)=1$.  There is some $M>0$
so that $f\circ\phi(t)>-M$ for all $t$.  

Let $B_R$ be the component of $f^{-1}[-R,+R]$ containing $1$.  If $R$ is
sufficiently large, then $B_R$ always separates $[g_1^\infty]$ from $[\phi]$ in
the sense that the two paths $g_1^\infty$ and $\phi$ are eventually in 
different components of the complement of $B_R$.  

We choose $R$ large enough so that $B_R$ separates $[g_1^\infty]$ from 
$[\phi]$ and so that $R$
is much larger than $M$ or $\cobf$.  
We may also choose $R$ so that $\phi$ crosses the
frontier of $B_R$ in an edge of $\Gamma(G,S)$.
Then the first group element
$h$ which $\phi$ passes through after leaving $B_R$ has $f(h)>R$.  Let
$g_2=hg_1^N$ where $N>\frac{99 R}{f(g)-\cobf}$, so that $f(g_2)>100R$. See
Figure  \ref{figure:lastbit}.
\begin{figure}[ht!]
\begin{center}
\input{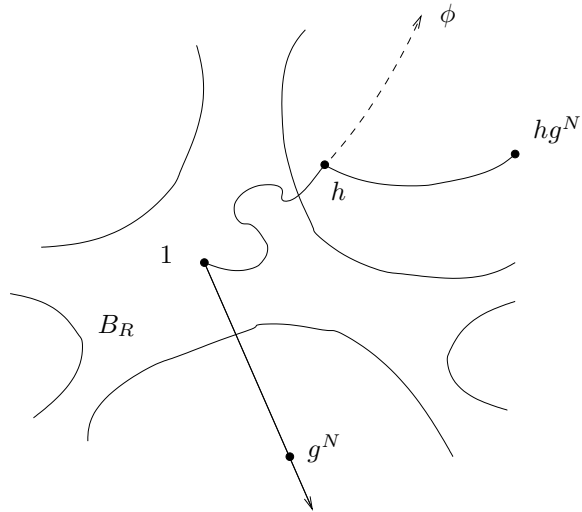}
\caption{$\phi$ may wiggle around a bit inside $B_R$ but the word representing
$g_2=hg^N$ is still ``coarsely monotone'' with respect to $f$.}
\label{figure:lastbit}
\end{center}
\end{figure}
We claim that $[g_2^\infty]$ is separated from $[g_1^\infty]$ by $B_R$.  We can
represent $g_2$ by a word $w=w_p g_1^N$ where $w_p$ is just the word traversed 
by the initial part of $\phi$.  Note that $f\circ\phi$ never decreases by more
than $2 R$ on this initial segment.
Let $\psi\co\R_+\to\Gamma(G,S)$ be the path
representing $[g_2^\infty]$ which traverses the infinite word $w^\infty$ at unit
speed starting at $1$.  If $\psi$ were to cross back over $B_R$ after getting to
$g_2$, we would have to have $f\circ\psi(t)< R$ for some $t>\mathrm{length}(w)$.
But since $f\circ\psi$ can decrease by no more than $2 R$, this is impossible.

Thus $[g_2^\infty]\neq [g_1^\infty]$ and both $f(g_1)$ and $f(g_2)$ are
positive.  The proof of the existence of $g_3$ is almost identical.
\end{proof}

\begin{theorem}\label{th:freegroup}
Let $G$ be a finitely presented group.  
If there is a bushy pseudocharacter on $G$,  
then $G$ contains a non-abelian free subgroup.  
\end{theorem}
\begin{proof}
Let $f\co G\to\R$ be a bushy pseudocharacter.  
To prove the theorem, it suffices to find $g$ and $g'$ in $G\setminus
f^{-1}(0)$ with disjoint fixed point sets in $E(f)$. 
If we can find such elements, then the Ping-Pong Lemma (see for instance
\cite[p467]{bridhaef:book})
 and  the dynamics described in Lemma
\ref{lemma:dynamics} ensure 
that high enough powers of these elements generate a
free group.

Let $g_1$, $g_2$, and $g_3$ be as in the proof of Lemma
\ref{lemma:groupelements}.  By taking powers we may assume that $f(g_i)>\cobf$
for $i\in\{1,2,3\}$.  If no two of these have disjoint fixed point sets,
then we must have $[g_3^\infty]=[g_1^\infty]=A$ and
$[g_2^{-\infty}]=[g_1^{-\infty}]=B$.  But then $g_2g_3(A)\neq A$ and
$g_2g_3(B)\neq B$, and so we may set $g=g_1$ and $g'=g_2g_3$.  These clearly
have disjoint fixed point sets.  Furthermore $f(g)=f(g_1)>0$ and $f(g')\geq
f(g_2)+f(g_3)-\cobf>0$.  
\end{proof}

\section{Quasi-actions on trees}\label{section:quasiaction}
In this section, we show that $G$ acts on a Gromov hyperbolic graph
quasi-isometric to a simplicial tree $\Gamma$, and that $E(f)$ embeds in the
ends of $\Gamma$.  If $f$ is not uniform, this implies that $G$ quasi-acts on
the bushy tree $\Gamma$ in the sense given in \cite{msw:quasiactI}:
\begin{definition}
A \emph{$(K,C)$--quasi-isometry} is a (not necessarily continuous) function
$q\co X\to Y$ between metric spaces so
that the following are true:
\begin{enumerate}
\item For all $x_1$, $x_2\in X$ 
\[d(x_1,x_2)/K-C\leq d(q(x_1),q(x_2))\leq Kd(x_1,x_2)+C.\]
\item The map $q$ is \emph{coarsely onto}, that is, every $y\in Y$ is 
distance at most $C$ from some point in $q(X)$. 
\end{enumerate}
\end{definition}
\begin{definition}
A \emph{$(K,C)$--quasi-action} of a group $G$ on a metric space $X$ is a map
$A\co G\times X\to X$, denoted $A(g,x)\mapsto gx$, so that 
the following hold:
\begin{enumerate}
\item For each $g$,
$A(g,-)\co G\to G$ is a $(K,C)$ quasi-isometry. 
\item For each $x\in X$ and $g$,
$h\in G$, we have $d(A(g,A(h,x)),A(gh,x))\leq C$.
(In other words, $d(g(hx),(gh)x)\leq C$.) 
\end{enumerate}
 We call a quasi-action \emph{cobounded} if
for every $x\in X$, the map $A(-,x)\co G\to X$ is coarsely onto.   
\end{definition}
\begin{definition}
Two quasi-actions $A_1\co G\times X\to X$ and $A_2\co G\times Y\to Y$ are called
\emph{quasi-conjugate} if there is a quasi-isometry $f\co X\to Y$ so that for
some $C\geq 0$ we have $d(f(A_1(g,x)), A_2(g,f(x))\leq C$ for all $x\in X$.  
The map $f$ is
called a \emph{quasi-conjugacy}.
\end{definition}
In contrast
to the quasi-actions discussed in \cite{msw:quasiactI}, the quasi-actions
on trees arising from pseudocharacters
are not in general quasi-conjugate to actions on trees.  We discuss ``exotic'' 
quasi-actions on trees further in Section
\ref{section:examples}. 

\subsection{Spaces Quasi-isometric to Trees}
It is helpful to develop a characterization of geodesic metric spaces
quasi-isometric to simplicial trees.  We will call a geodesic space a
\emph{quasi-tree} if it is quasi-isometric to some simplicial tree.  
One reason to be interested in quasi-trees is the following observation, which
was previously known to Kevin Whyte and probably to Gromov and others:
%\begin{proposition}
%Any quasi-action on a simplicial tree is quasi-conjugate to an isometric
% action on a quasi-tree.  Conversely, any isometric action on a quasi-tree is
%quasi-conjugate to a quasi-action on a simplicial tree.
%\end{proposition}
\begin{proposition}\label{prop:equiv}
Any quasi-action on a geodesic metric space $X$ is quasi-conjugate to
an action on some connected graph quasi-isometric to $X$.  Conversely, any isometric
action on a geodesic metric space quasi-isometric to $X$ is
quasi-conjugate to some quasi-action on $X$.
\end{proposition}
\begin{proof}[Sketch proof]
Suppose we have a $(K,C)$--quasi-action of a group $G$ on the space
$X$. Let $Y$ be a graph with
vertex set equal to $G\times X$, and connect $(g,x)$ to $(g',x')$ with an edge
(of length one)
whenever there is some $h\in G$ so that $d((hg)x,(hg')x')<2C$.  Define an action
of $G$ on the vertices of $Y$ by $g(h,x) = (gh,x)$.  Note that two vertices
connected by an edge will always be mapped to two vertices connected by an edge,
so this action extends to an isometric action on $Y$.  Let $f\co X\to Y$ be the
function $f(x)= (1,x)$.  It is not too hard to show that $f$ quasi-conjugates
the original quasi-action on $X$ to the action on $Y$.

Conversely, suppose that $X$ and $Y$ are quasi-isometric spaces, and suppose
$q\co Y\to X$ and $p\co X\to Y$  are $(K,C)$--quasi-isometries which are
$C$--quasi-inverses of one another (that is, $d(y,p(q(y)))$ and $d(x,q(p(x)))$
are bounded above by $C$ for all $y\in Y$ and $x\in X$).  Given an isometric
action of $G$
on $Y$, it is straightforward to check that $A\co G\times X\to X$
given by $A(g,x)= q(g(p(x)))$ is a $(K^2,KC+C)$--quasi-action.  
\end{proof}

In particular, any quasi-action on a simplicial tree is
quasi-conjugate to an
isometric
action on a quasi-tree and any isometric action on a quasi-tree is
quasi-conjugate to a quasi-action on a simplicial tree.

The following lemma is well known (see,
for example \cite[p401]{bridhaef:book}):
\begin{lemma}\label{lemma:quasigeodesic}
For all $K\geq 1$, $C\geq 0$, and $\delta\geq 0$ there is an $R(\delta,K,C)$ so that:

If $X$ is a $\delta$--hyperbolic metric space (e.g. a quasi-tree), 
$\gamma$ is a $(K,C)$--quasi-geodesic segment in $X$, and $\gamma'$ is a
geodesic segment with the same endpoints, then
the images of $\gamma'$ and $\gamma$ are Hausdorff distance less than $R$
from one another.  
\end{lemma}

\begin{theorem}\label{th:bottleneck}
Let $Y$ be a geodesic metric space. The following are equivalent:
\begin{enumerate}
\item  $Y$ is quasi-isometric to some simplicial tree $\Gamma$.
\item  {\rm(Bottleneck Property)}\qua There is some $\Delta>0$ so that 
for all $x$, $y$ in  $Y$ there is a midpoint $m=m(x,y)$ with 
$d(x,m)=d(y,m)=\frac{1}{2}d(x,y)$ and the property that
any path from $x$ to $y$ must pass within less than $\Delta$ of the point $m$.
\end{enumerate}
\end{theorem}
 
\begin{proof}\nl
\textbf{(1) $\Rightarrow$ (2)}\qua
Let $q\co Y\to \Gamma$ be a $(K,C)$--quasi-isometry, where $\Gamma$ is a 
simplicial tree.  Note that since $\Gamma$ is $0$--hyperbolic, and $Y$ is
quasi-isometric to $\Gamma$, $Y$ is $\delta$--hyperbolic for some $\delta$.
 Let $x$ and $y$ be two points of $Y$, joined by some geodesic segment 
$\gamma$.  Let $m$ be the midpoint of $\gamma$, and suppose that $\alpha$ is
some other path from $x$ to $y$.  

The image of a path under a $(K,C)$--quasi-isometry is a $C$--quasi-path.
In other words, though the path need not be continuous, it can make ``jumps'' of
length at most $C$.  
Therefore any point on the unique geodesic $\sigma$ from $q(x)$ to $q(y)$ in 
$\Gamma$ is no
more than $\frac{C}{2}$ from the image of $q\circ\alpha$.  
Furthermore $q\circ\gamma$ is a $(K,C)$--quasi-geodesic.  Thus by Lemma
\ref{lemma:quasigeodesic}, the distance from $q(m)$ to $\sigma$ is less than
$R=R(\delta,K,C)$.  

Let $p$ be the point on $\sigma$ closest to $q(m)$.  There is some point 
$z\in Y$ on
$\alpha$ so that $d(q(z),p)\leq\frac{C}{2}$.  Since $d(p,q(m))<
R$ we have \[d(q(z),q(m))< C/2+R\] which implies 
\[d(z,m)< K(C/2+R)+C.\]
In other words, the path $\alpha$ must pass within
$K\bigl(\frac{C}{2}+R\bigr)+C$ of the
point $m$, so we may set $\Delta=K\bigl(\frac{C}{2}+R\bigr)+C$.

\textbf{(2) $\Rightarrow$ (1)}\qua
Given a geodesic metric space with the Bottleneck 
Property, we inductively construct a
simplicial tree and a quasi-isometry from it to $Y$.  At each stage of the
construction we have a map $\beta_k\co \Gamma_k\to Y$ where $\Gamma_k$ is a tree
of diameter $2k$.  We let $V_k$
be the image of the vertices of $\Gamma_k$ under $\beta_k$, and let
$N_k$ be a large
neighborhood of $V_k$.  We refer to the set of path 
components of the complement of
$N_k$ as $\mathcal{C}_{k+1}$.  Each element of $\mathcal{C}_{k+1}$ gives rise to a vertex
in $\Gamma_{k+1}\setminus \Gamma_k$.

\textbf{Step 0}\qua
Let $R=20 \Delta$.  Pick some base point $*\in Y$.  We set $V_0=\{*\}$, and
$\Gamma_0$ equal to a single point $p_*$.  For each $i$ there will be a natural
identification of $V_i$ with the vertices of $\Gamma_i$.  We define 
$\beta_0\co\Gamma_0\to Y$ so that the image of $\beta_0$ is $*$.  We define
$\mathcal{C}_0=\{Y\}$.

\textbf{Step k}\qua
Let $k\geq 1$, and define $d_i\co Y\to \R$ by $d_k(x)=d(x,V_{k-1})$.
Let $N_{k-1}=\{x\in Y\ |\ d_i(x)<R\}$
and let $\mathcal{C}_k$ be the set of path components of
$d_i^{-1}\left([R,\infty)\right)$.
If $C\in\mathcal{C}_k$, let $\mathrm{Front}(C)=\{x\in C\ |\
d_i(x)=R\}$.  Because $d_i$ is continuous and $Y$ is path-connected, 
$\mathrm{Front}(C)$ is nonempty.
We pick some $v_C\in \mathrm{Front}(C)$ for each $C$ in
$\mathcal{C}_k$, and let  $U_k=\{v_C\ |\ C\in\mathcal{C}_k\}$.  

The new set $V_k$ is equal to $V_{k-1}\cup U_k$. 
 To construct $\Gamma_k$, we add one
vertex $p_v$ to $\Gamma_{k-1}$ for each $v\in U_k$.  The point $v$ is contained
in exactly one element of $\mathcal{C}_{k-1}$, and this element
contains exactly one
element $w$ of $V_{k-1}$.  We connect the new vertex $p_v$ to the old vertex
$p_w$ by a single edge.  The map $\beta_k$ is defined to be equal to $\beta_{k-1}$ on
$\Gamma_{k-1}$, and is extended to map a new edge between $p_w$ and $p_v$ to a
geodesic segment in $Y$ joining $w$ to $v$.  This completes Step $k$.  
See Figure \ref{figure:buildgamma} for an example of what this might look like.
\begin{figure}[ht!]
\begin{center}
\input{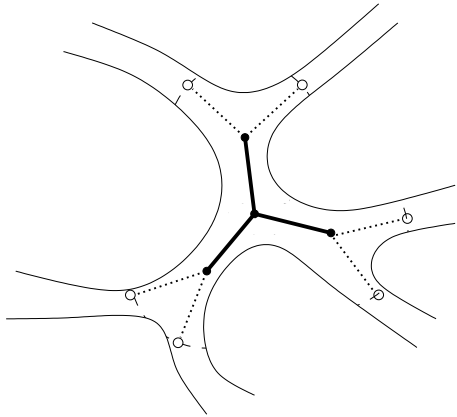}
\caption{A stage in the construction of $\Gamma$. 
The image of $\Gamma_{k-1}$ is a thick tripod.  The elements of
$U_k$ are open circles, and the new edges of $\Gamma_k$ are mapped to the dotted
segments.}
\label{figure:buildgamma}
\end{center}
\end{figure}

After the induction is completed, we have a map 
$\beta\co\bigcup_k\Gamma_k=\Gamma\to Y$ which is defined to be equal to $\beta_k$ on
each $\Gamma_k$.  The image of the $0$--skeleton of $\Gamma$ in $Y$ is
$V=\bigcup_k V_k$.  
We will show that $\beta\co\Gamma\to Y$ is a quasi-isometry by showing that
$\beta$
restricted to the $0$--skeleton of $\Gamma$ is a quasi-isometry.  We will use
the following lemma:
\begin{lemma}\label{lemma:sixdelta}
Let $v\in U_i=V_i\setminus V_{i-1}$, 
and suppose that $p_v$ is connected to $p_w\in \Gamma_{i-1}$ by
an edge.  Then the following assertions hold:
\begin{enumerate}
\item $R\leq d(v,w)\leq R+6\Delta$
\item If $v\in C\in \mathcal{C}_i$ and 
$p\in \mathrm{Front}(C)$, then 
$d(v,p)\leq 6\Delta$.
\end{enumerate}
\end{lemma}
\begin{proof}
We prove both assertions simultaneously by induction.

\textbf{Step 1}\qua  If $i=1$, then $V_{i-1}=V_0$ is a single point,
so assertion
(1) holds with $d(v,w)=R$.

We prove assertion (2) by way of contradiction.  Namely, suppose that
$d(v,p)>6\Delta$.  There is a point $m$ equidistant from $v$ and $p$ so that any
path from $v$ to $p$ passes within less than $\Delta$ of $m$.  As there is a
path $\beta$ in $C$ connecting $v$ to $p$, we must have
 $d(m,C)=\inf_{c\in C}\{d(m,c)\}<\Delta$.  The situation is shown in Figure
\ref{figure:step0contradict}.
\begin{figure}[ht!]
\begin{center}
\input{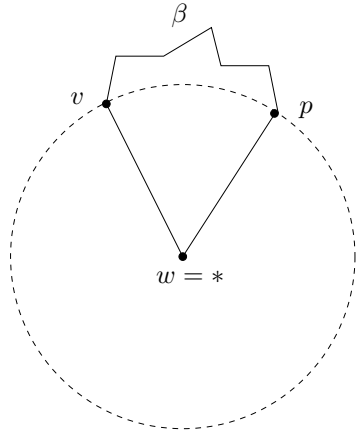}
\caption{Assertion (2) Step 1.  Where is the midpoint between $p$ and $v$?}
\label{figure:step0contradict}
\end{center}
\end{figure}

On the other hand, there is another path from $v$ to $p$ consisting of a
geodesic segment $[v,w]$ from $v$ to $w$ and another geodesic segment $[w,p]$
from $w$ back to $p$.  By the Bottleneck property, $m$ must lie inside a
$\Delta$--neighborhood of this path.  In other words, there is some point $z$ 
on the path $[v,w]\cup[w,p]$ so that $d(z,m)<\Delta$.  Since $d(v,p)>6\Delta$,
it follows that
$d(v,m)$ and $d(m,p)$ are both strictly greater than $3\Delta$, and thus
$d(z,C)=d(z,\{v,p\})>2\Delta$ by the triangle inequality.  Since
$d(m,z)<\Delta$, we have $d(m,C)>2\Delta-\Delta= \Delta$, a contradiction.

\textbf{Step i}\qua  We again prove assertion (1) first, now assuming that 
$v\in U_i$ is in the same $D\in\mathcal{C}_{i-1}$ as $w\in V_{i-1}$.  Because
$d(v,V_{i-1})=R$, for any $\epsilon>0$ there is some $w'\in V_{i-1}$ with
$R\leq d(v,w')\leq R+\epsilon$.  If for all $\epsilon>0$ we can choose $w'=w$,
then $d(v,w)=R$ and we are done, so assume that $w'\neq w$.  Then $w'$ is not
contained in $D$, so a geodesic path from $v$ to $w'$ must pass through some
point $d\in\mathrm{Front}(D)$ (see Figure \ref{figure:step_i_assertion_1}).
\begin{figure}[ht!]
\begin{center}
\input{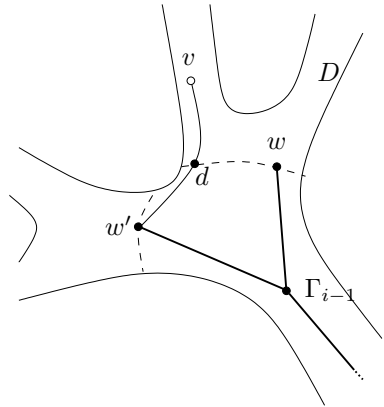}
\caption{Assertion (1) Step i}
\label{figure:step_i_assertion_1}
\end{center}
\end{figure}
By the induction hypothesis (2),
$d(d,w)\leq6\Delta$, so $d(v,w)\leq R+\epsilon+6\Delta$.  Letting $\epsilon$
tend to zero, we obtain assertion (1). 

To prove assertion (2) we again argue by way of contradiction.  Let
$p\in\mathrm{Front}(C)$ be such that $d(p,v)>6\Delta$.  We see as before that
the midpoint $m$ provided by the Bottleneck Property must satisfy
$d(m,C)<\Delta$.  

To obtain the contradictory inequality in this case requires some extra 
maneuvers.  Let $\epsilon > 0$.  Then we may find $w_1$ and $w_2$ in
$V_{i-1}$ so that
$d(v,w_1)\leq R+\epsilon$ and $d(p,w_2)\leq R+\epsilon$.  These points $w_1$ and
$w_2$ are connected by a path $\sigma$ in $\beta(\Gamma_{i-1})$ which is the image
of a geodesic path in the tree $\Gamma_{i-1}$.  Together with geodesics
$[v,w_1]$ and $[w_2,p]$, this path $\sigma$ gives a path between $v$ and $p$
(see Figure \ref{figure:step_i_assertion_2a}).
\begin{figure}[ht!]
\begin{center}
\input{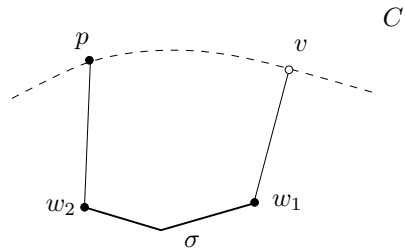}
\caption{Assertion (2) Step i}
\label{figure:step_i_assertion_2a}
\end{center}
\end{figure}
The midpoint $m$ must lie in a $\Delta$--neighborhood of this path.  Let $z$ be
a point on the path $[v,w_1]\cup\sigma\cup[w_2,p]$ which is less than $\Delta$
from $m$.  We claim first that $z$ cannot lie on $\sigma$.  Certainly $z$ cannot
be an element of $V_{i-1}$, as $d(V_{i-1},C)=R=20\Delta$.  Suppose then that $z$
is in the interior of an edge of $\sigma$.  Using the triangle inequality and
the induction hypothesis (each such edge must have length between $R$ and
$R+6\Delta$), we see that $d(z,C)\geq 7\Delta$, so we would have $d(m,C)\geq
6\Delta$, a contradiction.

The only remaining possibility is that $z$ lies on one of the geodesic segments
$[v,w_1]$ or $[w_2,p]$.  We may suppose that $z$ lies on $[v,w_1]$.  Assume that
$q$ is an arbitrary point in $C$.  We will argue from the triangles shown in
Figure \ref{figure:step_i_assertion_2b}.
\begin{figure}[ht!]
\begin{center}
\input{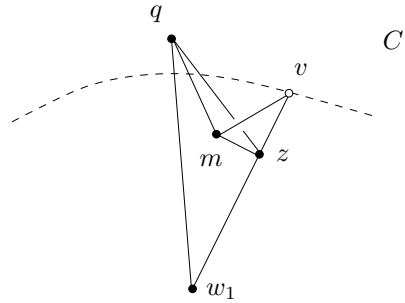}
\caption{Assertion (2) Step i (continued)}
\label{figure:step_i_assertion_2b}
\end{center}
\end{figure}
By assumption, $d(z,m)<\Delta$.  Since $d(v,m)>3\Delta$, $\triangle vmz$ gives 
$d(z,v)>2\Delta$.  We also have $d(v,w_1)\leq R+\epsilon$, and so since $z$ lies
on a geodesic from $w_1$ to $v$
we get that $d(z,w_1)< R+\epsilon-2\Delta$.  Since $q\in C$, we have
$d(w_1,q)\geq R$, and so $\triangle qzw_1$ gives $d(z,q)>2\Delta-\epsilon$. 
 Finally $\triangle mqz$ gives $d(m,q)>\Delta-\epsilon$. 
Letting $\epsilon$ tend to zero gives
$d(m,C)\geq \Delta$.  This contradiction establishes Assertion (2).
\end{proof}

%%%%% coarsely onto
That $\beta$ is coarsely onto follows easily from Lemma
\ref{lemma:sixdelta}.  Indeed, suppose 
that $x\in Y$ is not contained in an $R$--neighborhood of $V\subset
\beta(\Gamma)$. 
Because $Y$ is geodesic, we may find an $x$ which is
distance exactly $R$ from $V$.  There is then some $i$ and some $v\in V_i$ so 
that $d(x,v)< R+\Delta$.  Of course, $x$ must lie 
in some component
of the complement of $N_i = \{x\ |\ d(x,V_i)<R\}$, in other words there is some
$C\in\mathcal{C}_{i+1}$ with $x\in C$.  
Let $w_C$ be the element of $U_{i+1}$ corresponding to this
component.  A geodesic path from $x$ to $v$ must pass through
$\mathrm{Front}(C)$ at some point $p$.  
Since $d(x,v)< R+\Delta$ and $d(v,C)\geq R$, we have $d(p,x)<\Delta$. 
But by Lemma \ref{lemma:sixdelta},
$d(p,w_C)< 6 \Delta$, so $d(x,w_C)< 7\Delta < R$, contradicting our choice of
$x$.
%%%%% 

%%%%% upper bound
The images of the edges of $\Gamma$ allow us to get an upper bound on
$d(\beta(x),\beta(y))$ where $x$ and $y$ are vertices of $\Gamma$.  By Lemma
\ref{lemma:sixdelta} the image of each such edge is a geodesic of length less
than or equal to $R+6\Delta$.  Thus 
\[d(\beta(x),\beta(y))\leq (R+6\Delta)d(x,y)=26\Delta d(x,y).\]

%%%%% lower bound
Let $x$ and $y$ be vertices of $\Gamma$.  These are joined in $\Gamma$ by a
unique geodesic $\sigma$.  For any $p\in\Gamma$ we can define $D(p)$ to be the
minimum $i$ so that $p\in\Gamma_i$.  If $p$ on $\sigma$ minimizes $D$ and is not
an endpoint of $\sigma$, we refer to $p$ as the \emph{turnaround vertex}.  Note
that $\sigma$ contains at most one turnaround vertex, as $\Gamma$ is simply
connected.  
\begin{lemma}\label{lemma:lowerbound}
If $z$ is a vertex on $\sigma$ but is not a turnaround vertex, then any geodesic
from $\beta(x)$ to $\beta(y)$ passes within $6\Delta$ of $\beta(z)$.
\end{lemma}
\begin{proof}
Without loss of generality, we may assume $D(x)>D(z)$ and $D$ is non-increasing
from $x$ to $z$. Otherwise we may switch $x$ and $y$ to ensure this is the
case.  
Let $C\in\mathcal{C}_{D(z)}$ be the component of 
$Y\setminus N_{D(z)-1}$ containing $z$.  We claim that $\beta(x)\in C$ but 
$\beta(y)\notin C$.

We see that $\beta(x)\in C$ by induction on $d(x,z)$.  Let $z=z_0,z_1,\ldots,z_N=x$
be the sequence of vertices on $\sigma$ joining $z$ to $x$.  Let
$C=C_0,C_1,\ldots,C_N$ be the corresponding complementary components
$C_i\in\mathcal{C}_{D(z)+i}$ from the construction of $\Gamma$.  For $i\geq 1$
we have $C_i\subset C_{i-1}$ by construction, so $\beta(x)\in C$.  

If $\beta(y)$ were in $C$, then there would be a path from $y$ back to $z$ on which $D$ was non-increasing, and thus $z$ would be a turnaround vertex.

Since $\beta(x)\in C$ but $\beta(y)\notin C$, any geodesic from $\beta(x)$ to $\beta(y)$ must
pass through $\mathrm{Front}(C)$, and so by Lemma \ref{lemma:sixdelta}, the
geodesic must pass within $6\Delta$ of $z$.
\end{proof}

By Lemma \ref{lemma:lowerbound} a geodesic from $\beta(x)$ to $\beta(y)$ must pass
within $6\Delta$ of each $\beta(z)$ where $z$ is a vertex of $\Gamma$ between $x$
and $y$ which is not the turnaround vertex.  Images under $\beta$ of successive
vertices are at least $R$ apart, so by picking points on the geodesic with
$6\Delta$ of the images of the vertices (except for the turnaround vertex) we
see that \[d(\beta(x),\beta(y))\geq (R-12\Delta)(d(x,y)-2) = 8\Delta
d(x,y)-16\Delta.\]
Combining this result with the previously obtained upper bound we get
\[8\Delta d(x,y)-16\Delta\leq d(\beta(x),\beta(y))\leq 26\Delta d(x,y).\]
In particular, $\beta$ is a quasi-isometric embedding which is $R$--almost onto, so
it is a quasi-isometry.
\end{proof}

\subsection{Pseudocharacters and Quasi-actions} 
Just as a homomorphism $\chi\co G\to \R$ gives rise to a $G$--action on $\R$
via $g(x)=\chi(g)+x$ for $g\in G$ and $x\in\R$, a pseudocharacter $f\co G\to \R$
gives rise to a $(1,\cobf)$--quasi-action of $G$ on $\R$ via $g(x)=f(g)+x$.  
Roughly, in
this section we attempt to ``lift'' this quasi-action to a quasi-action on $T$,
the tree defined in Section \ref{section:topology}. It is not immediately clear whether the quasi-action should lift,
because the image of a vertex space (in $\widetilde{K}$) after 
action by a group
element might intersect infinitely many vertex spaces.  The point of the
construction in this section is that we may ``collapse'' enough of $T$ to get a
quasi-action, and still be left with a complicated enough tree so that $E(f)$
embeds in its Gromov boundary.

Recall the definition of the tree $T$.    
We pick an 
(unambiguous) triangular generating set $S$. We then scale $f$
so that $f(G)$ misses $\Z+\frac{1}{2}$ and so that 
$f$ changes by at most $\frac{1}{4}$ over each edge.  
We then build a tree with vertex
set in one-to-one correspondence with the components of $\widetilde{K}\setminus
f^{-1}\bigl(\Zhalf\bigr)$.  The edges correspond to components of
$f^{-1}\bigl(\Zhalf\bigr)$, each of which is some possibly infinite track which separates
$\widetilde{K}$ into two components.

\begin{definition}\label{definition:X} 
We now define a graph 
$X$, which we will later show is a quasi-tree.  
Let $V$ be the set of components of $f^{-1}\bigl(\Zhalf\bigr)$. Then $V$ is in 
one-to-one correspondence with the set of edges of $T$. 
Let $X$ be the simplicial graph with vertex set equal to $G\times V$ and the
following edge condition:  Two distinct vertices $(g,\tau)$ and $(g', \tau')$ 
are to be
connected by an edge if there is some $h$ so that $hg(\tau)$ and $hg'(\tau')$ 
are contained in the same
 connected component of $f^{-1}\bigl[n-\frac{3}{2}, n+\frac{1}{2}\bigr]$ for
some $n\in\Z$.  We endow the zero-skeleton $X^0$ with a $G$--action by setting
$g(g_0,\tau_0)=(gg_0,\tau_0)$.  Since this action respects the edge condition on
pairs of vertices, it extends to an action on $X$.
\end{definition}
\begin{remark}
The relationship between $X$ and $T$ is actually somewhat unclear.
For every $x\in T$ we choose $e(x)$ to be an arbitrary edge adjacent
to $x$.
If $\varpi\co T\to X$ is given by $\varpi(x)=(1,e(x))$, then $\varpi$
is coarsely surjective and coarsely Lipschitz (in fact,
$d(\varpi(x),\varpi(y))\leq d(x,y)+1$).  If $\varpi$ were a
quasi-isometry, then
Theorem \ref{th:pseudoquasibushy} would follow immediately: $G$ would
quasi-act on the tree $T$.  Usually, though, $\varpi$ is not a
quasi-isometry;  the preimages of bounded sets do not even need to be
bounded.  
\end{remark}
\begin{lemma}\label{lemma:connected}
$X$ is connected.
\end{lemma}
\begin{proof}
For every edge $e$ of $T$, there is a vertex $(1,\tau_e)$ where $\tau_e$ is the
component of $f^{-1}\bigl(\Zhalf\bigr)$ corresponding to $e$.
If $e_1$ and $e_2$ are adjacent edges of $T$, then $(1,\tau_{e_1})$ and 
$(1,\tau_{e_2})$ are
certainly connected by an edge in $X$.  Thus the vertices $\{1\}\times V$ are
all in the same connected component of $X$.  

Let $(g,\tau)$ be some
vertex of $X$.  As $G$ acts by isomorphisms of the complex $\widetilde{K}$, $g\tau$
is, like $\tau$, some track in $\widetilde{K}$.  For any point
$x\in\widetilde{K}$,
$|f(gx)-f(g)-f(x)|\leq \cobf$, since $f$ on $\widetilde{K}$ is obtained from $f$ on
$G$ by affinely extending over each cell.
In particular, since $f(\tau)$ is a point, $f(g\tau)$ has diameter less
than or equal to $2\cobf$ in $\R$.  Since by assumption $\cobf$ is much less
than one, $f(g\tau)$ is contained in the interval 
$\bigl[n-\frac{3}{2}, n+\frac{1}{2}\bigr]$ for
some $n\in\Z$.  Since $g\tau$ is connected, it is therefore contained in a
connected component of $f^{-1}\bigl[n-\frac{3}{2}, n+\frac{1}{2}\bigr]$.  The boundary of
this set contains at least one component $\tau'$ of $f^{-1}(\Z)$, and so
$(g,\tau)$ is connected to $(1,\tau')$ by an edge of $X$.  Thus all the vertices
of $X$ are contained in the same connected component, and $X$ is connected.
\end{proof}

$G$ clearly acts simplicially on $X$. 
If we regard $X$ as a path metric space with each edge having length $1$, then
$G$ acts isometrically on $X$. 
\newcommand{\gtau}[1]{\ensuremath{(g_{#1},\tau_{#1})}}
\begin{proposition}
$G$ acts coboundedly on $X$.
\end{proposition} 
\begin{proof}
Let $\gtau{0}$ be a vertex of $X$.  We will show that every other vertex of $X$
is distance at most $1$ from the orbit of $\gtau{0}$.  Let $\gtau{1}$ be another
vertex of $X$.  For $i\in \{0,1\}$, let $e_i$ be an edge which intersects
$g_i\tau_i$, and let $h_i\in G$ be an endpoint of $e_i$.  Then
$h_i^{-1}g_i\tau_i$ is a track which passes through an edge adjacent to $1$.
Thus $\sup|f(h_i^{-1}g_i\tau_i)|\leq \epsilon_f+\cobf<2\epsilon_f<\frac{1}{2}$.
Since both tracks pass through edges adjacent to $1$, a single component of
$f^{-1}\bigl[-\frac{3}{2},\frac{1}{2}\bigr]$ contains both $h_0^{-1}g_0\tau_0$ and
$h_1^{-1}g_1\tau_1$.  According to Definition \ref{definition:X}, this means that
$d((h_0^{-1}g_0,\tau_0),(h_1^{-1}g_1,\tau_1))\leq 1$.  Since $G$ acts on $X$ by
isometries, $d(h_1 h_0^{-1}\gtau{0},\gtau{1})\leq 1$ and the proposition is proved.
\end{proof}

\begin{lemma}\label{lemma:intersect}
Suppose $\gtau{1}$ and $\gtau{2}\in X$ and suppose that $g_1\tau_1\cap
g_2\tau_2$ is nonempty.  Then 
%either $\gtau{1}=\gtau{2}$ or there is an edge
%connecting $\gtau{1}$ to $\gtau{2}$.  In other words,
$d(\gtau{1},\gtau{2})\leq 1$.  
\end{lemma}
\begin{proof}
Suppose that $\gtau{1}$ and $\gtau{2}$ are distinct, and let $h$ be a vertex of
some $2$--cell in $\widetilde{K}$ through which $g_1\tau_1$ and $g_2\tau_2$ both
pass.  Both $h^{-1}g_1\tau_1$ and $h^{-1}g_2\tau_2$ pass through a $2$--cell
adjacent to $1$, and so 
$\sup |f(h^{-1}g_1\tau_1\cup h^{-1}g_2\tau_2)|<\epsilon_f<\frac{1}{4}$.  Thus
$h^{-1}g_1\tau_1\cup h^{-1}g_2\tau_2$ is contained in a single component of
$f^{-1}\bigl[-\frac{3}{2},\frac{1}{2}\bigr]$, and
$d(\gtau{1},\gtau{2})=d(h^{-1}\gtau{1},h^{-1}\gtau{2})=1$.
\end{proof}

\begin{lemma}\label{lemma:Tbottle}
Let $(g_a,\tau_a)$, $(g_b,\tau_b)$, $(g_c,\tau_c)\in X$ be such that $g_b\tau_b$
separates $g_a\tau_a$ from $g_c\tau_c$ in $\widetilde{K}$.  Then any path from
$(g_a,\tau_a)$ to $(g_c,\tau_c)$ passes within $2$ of $(g_b,\tau_b)$.
\end{lemma}
\begin{proof}
Let $\gtau{a}=\gtau{0},\gtau{1},\ldots,\gtau{n}=\gtau{c}$ be the vertices of a
path in $X$ connecting \gtau{a}\ to \gtau{c}.  If $g_k\tau_k$ intersects
$g_b\tau_b$ for any $k$ then we have $d(\gtau{k},\gtau{b})\leq 1$ by Lemma
\ref{lemma:intersect}.  Thus we may assume that $g_k\tau_k$ is disjoint from
$g_b\tau_b$ for all $k$.  Since $g_b\tau_b$ separates $g_0\tau_0$ from
$g_n\tau_n$ there is some $k$ for which $g_k\tau_k$ and $g_{k+1}\tau_{k+1}$ are
separated by $g_b\tau_b$.  Since \gtau{k}\ is connected to \gtau{k+1}\ by an
edge of $X$, there is some $h\in G$, some $n\in\Z$, and some connected component
$B$ of $f^{-1}\bigl[n-\frac{3}{2},n+\frac{1}{2}\bigr]$ so that $hg_k\tau_k\cup
hg_{k+1}\tau_{k+1}\subset B$.  Since $B$ is path-connected, there is some 
path $\gamma\co[0,1]\to B$ with
$\gamma(0)\in hg_k\tau_k$ and $\gamma(1)\in hg_{k+1}\tau_k$.  Of course this
path must cross $hg_b\tau_b$, so $hg_b\tau_b\cap B$ is nonempty.  If
$hg_b\tau_b$ were contained in $B$, we would have $hg_b\tau_b\cup
hg_k\tau_k\subset B$, and so \gtau{b}\ and \gtau{k}\ would be connected by an
edge, implying $d(\gtau{b},\gtau{k})=1$.  If on the other hand $hg_b\tau_b$ is
not completely contained in $B$, then it must intersect some boundary component
$\tau$ of $B$.  Since $\tau=1\cdot\tau$ and $hg_b\tau_b$ intersect, we have
$d((1,\tau),(hg_b,\tau_b))\leq 1$ which implies $d((h^{-1},\tau),\gtau{b})\leq
1$.  Since $\tau\cup hg_k\tau_k\subset B$, we also have
$d((h^{-1},\tau),\gtau{k})=1$.  By the triangle inequality, 
$d(\gtau{k},\gtau{b})\leq 2$, establishing the lemma.  
\end{proof}

\begin{theorem}\label{th:xbottle}
The space $X$ satisfies the Bottleneck Property of Theorem \ref{th:bottleneck}
for $\Delta = 10$.
\end{theorem}
\begin{proof}
Let $x$, $y\in X$, and let $m$ be the midpoint of some geodesic segment $\gamma$
joining $x$ to $y$.  We may assume that $d(x,y)>20$, otherwise the Bottleneck
Property is satisfied trivially for $\Delta = 10$.  

The proof of Lemma \ref{lemma:connected} shows that any vertex of $X$ is
distance at most $1$ from some vertex of the form $(1,\tau)$.  As any point of
$X$ is distance at most $\frac{1}{2}$ from some vertex, there exist 
$(1,\tau)$ and $(1,\tau')$ so that
$d(x,(1,\tau))\leq \frac{3}{2}$ and $d(y,(1,\tau'))\leq \frac{3}{2}$.

Let $e$, $e'$ be the edges in $T$ corresponding to $\tau$ and $\tau'$, and let
$m_e$, $m_{e'}$ be the midpoints of these edges.  These points are connected by
a unique geodesic passing through some sequence of edges
$e=e_0,e_1,\ldots,e_n=e'$ in $T$.  If $\tau_i$ is the track in $\widetilde{K}$
associated to $e_i$, note that $d((1,\tau_i),(1,\tau_j))\leq |i-j|$.  In
particular, $d((1,\tau_i),(1,\tau_{i+1}))=1$ for all $i$.  Thus the sequence of
edges in $T$ defines a path in $X$ leading from $(1,\tau)=(1,\tau_0)$ to
$(1,\tau')=(1,\tau_n)$.  We extend this path with geodesic segments of length
less than or equal to $\frac{3}{2}$ to obtain a path from $x$ to $y$, most of
whose vertices lie in $1\times V$.  This path clearly contains some point $z$ so
that $\min\{d(x,z),d(y,z)\}\geq \frac{d(x,y)}{2}>10$.  
Thus there is some vertex $(1,\tau_k)$ such that
\[\min\{d(x,(1,\tau_k)),d(y,(1,\tau_k))\}\geq\frac{d(x,y)}{2}-\frac{1}{2}>\frac{19}{2}.\]
Since $(1,\tau_k)$ is far from $x$ and $y$, $k$ is not equal to $0$ or $n$, and
so $\tau_k$ separates $\tau$ from $\tau'$ in $\widetilde{K}$.  Thus by Lemma
\ref{lemma:Tbottle} any path from $(1,\tau)$ to $(1,\tau')$ must pass within
$2$ of $(1,\tau_k)$.  

Let $\sigma$ be any path from $x$ to $y$.  This path can be extended by adding
segments of length at most $\frac{3}{2}$ at each end to give a path
$\overline{\sigma}$ from $(1,\tau)$ to $(1,\tau')$.  By the previous paragraph,
$\overline{\sigma}$ must pass within $2$ of $(1,\tau_k)$.  Since the appended
segments are very far (at least $8$) from $(1,\tau_k)$, this means that $\sigma$
must pass within $2$ of $(1,\tau_k)$.  

By the same argument, the geodesic $\gamma$ passes within $2$ of $(1,\tau_k)$.
Let $z$ be a point on $\gamma$ which is within $2$ of $(1,\tau_k)$.  By the
triangle inequality,
 $\min\{d(x,z),d(y,z)\}\geq \frac{d(x,y)}{2}-\frac{5}{2}$.
Thus $d(z,m)\leq \frac{5}{2}$, where $m$ is the midpoint of $\gamma$.  Thus
$d(\sigma,m)\leq \frac{11}{2}< 10$, establishing the theorem. 
\end{proof}

\begin{corollary}\label{cor:qi}
The space $X$ is quasi-isometric to a simplicial tree $\Gamma$.
Since $G$ acts coboundedly on $X$, $G$ quasi-acts coboundedly on $\Gamma$.
\end{corollary}

%Of course, this is somewhat less impressive if $\Gamma$ is two-ended. 
%(Note that $\Gamma$
%must have at least two ends, because there is a map from $X$ onto $\R$ which 
%cannot stretch distances very much --
% $(g,\tau)\mapsto f(g)+f(\tau)$ gives such a map.)

For the next lemma, we need the language of Gromov products in metric spaces.
Recall that if $(M,m)$ is a pointed metric space, and $x$, $y\in M$, then the
Gromov product of $x$ with $y$ is defined to be
$(x.y)=\frac{1}{2}(d(m,x)+d(m,y)-d(x,y))$. If $M$ is Gromov hyperbolic, we say
that a sequence $\{x_i\}$ of points in $M$ \emph{converges at infinity} if
$\lim_{i,j\to\infty}(x_i.x_j)=\infty$.  One can then define $\partial M$,
 the Gromov boundary of $M$, to be the set of sequences converging at infinity
modulo the equivalence relation: $\{x_i\}\sim\{y_i\}$ if
$\lim_{i,j\to\infty}(x_i.y_j)=\infty$.  For more detail, see
 \cite[Chapter~III.H]{bridhaef:book} and \cite{gromov:wordhyperbolic}.  
\begin{lemma}\label{lemma:injective}
There is an injective map from $E(f)$ to $\partial X$.
\end{lemma}
\begin{proof}
Fix some base point $(1,\nu)\in X$.  All Gromov products in $X$ 
will be taken with
respect to this base point.  Let $[\phi]\in E(f)$.  The path $\phi$
passes through some sequence of tracks in
$f^{-1}\bigl(\Z+\frac{1}{2}\bigr)$.  We choose a
subsequence $\{\tau_i\}$ of these tracks so that for each $i$ the track $\tau_i$
separates $\nu$ from $\tau_{i+1}$.  Recall that these $\tau_i$ may be identified
with edges of the tree $T$ from Section \ref{section:topology}.  
Choosing them the
way we have ensures that they all lie on a geodesic ray in $T$.  It should be
clear that this geodesic ray limits on $i([\phi])$ where $i$ is the map from
Proposition \ref{prop:tree}.

Because $[\phi]\in E(f)$, we necessarily
have $\lim_{i\to\infty}f(\tau_i) = \pm \infty$.  We claim that the sequence
$\{x_i\}=\{(1,\tau_i)\}$ converges at infinity.  To show this we must show that
$(x_i.x_j)\to\infty$.  Let $0<i\leq j$.  Then 
\[(x_i.x_j)=\frac{1}{2}\bigl[d(x_i,(1,\nu))+d(x_j,(1,\nu))-d(x_i,x_j)\bigr].\]
If $i=j$ then clearly $(x_i.x_j)=d(x_i,(1,\nu))$.  Otherwise, $\tau_i$ is
between $\nu$ and $\tau_j$, so any path (in particular a geodesic) from
$(1,\nu)$ to $x_j$ must pass within distance $2$ of $x_i$, by Lemma
\ref{lemma:Tbottle}.  Thus we have 
\[d(x_j,(1,\nu))\geq d(x_i,x_j)+d(x_i,(1,\nu))-4\]
which implies
\[(x_i.x_j)\geq \frac{1}{2}\bigl(2 d(x_i,(1,\nu))-4\bigr)=d(x_i,(1,\nu))-2.\]
Since $d(x_i,(1,\nu))\geq \frac{1}{2}|f(\tau_i)-f(\nu)|\to\infty$ as
$i\to\infty$, so does $(x_i.x_j)$ and so $\{x_i\}$ converges to some point in
$\partial X$.  An almost identical calculation shows that this point does not
depend on the
sequence of tracks we chose.  
We define a map $\overline{i}\co E(f)\to \partial X$ by defining
$\overline{i}([\phi])$ to be the point in $\partial X$ to which this sequence 
converges.

We next claim that $\overline{i}$ is injective.  Suppose that $[\phi]$ and $[\phi']$
are distinct elements of $E(f)$.  By Proposition \ref{prop:tree} we may identify
$[\phi]$ and $[\phi']$ with elements of $\partial T$.  There is a unique
bi-infinite sequence of edges of $T$ (ie a geodesic) joining $[\phi]$ to 
$[\phi']$, so that one of every triple of edges separates the remaining two from
one another.  We will abuse notation slightly by identifying an edge of $T$
with the track associated to it.  Let $\omega$ be the edge adjacent to the 
geodesic between $[\phi]$ and $[\phi']$ which is closest to $\nu$.
\begin{figure}[ht!]
\begin{center}
\input{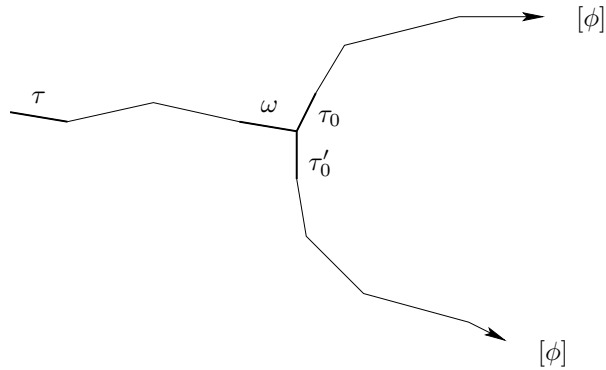}
\caption{Arrangement of edges in $T$}
\label{figure:inject}
\end{center}
\end{figure}
We may choose the sequences $x_i=(1,\tau_i)$ and $x_i'=(1,\tau_i')$ 
representing $\overline{i}([\phi])$ and $\overline{i}([\phi'])$ so that $\tau_0$ and
$\tau_0'$ are as shown in Figure \ref{figure:inject}.  Let $i$, $j\geq 1$.  Then
$\tau_i$ is separated by $\tau_0$ from $\nu$ and $\tau_j'$ and $\tau_j'$ is
separated by $\tau_0'$ from $\tau_i$ and $\nu$.  Transporting this arrangement
into $X$ we get Figure \ref{figure:injecttriang}.
\begin{figure}[ht!]
\begin{center}
\input{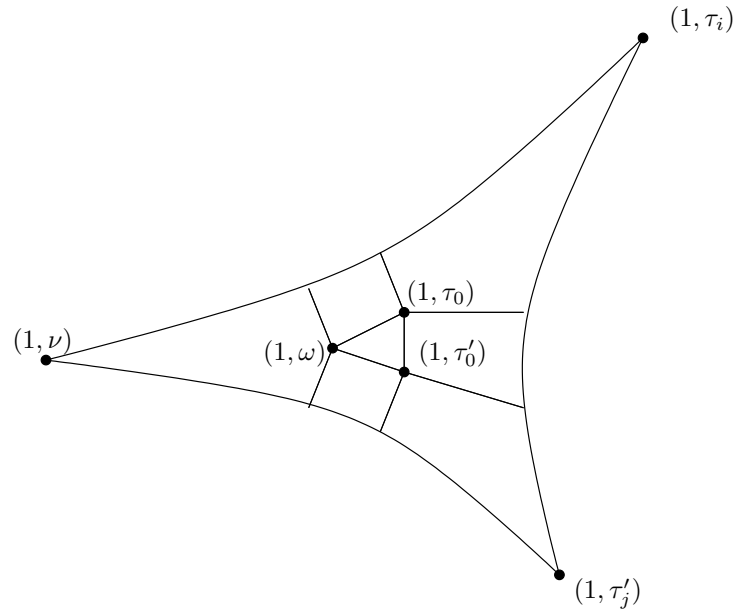}
\caption{Computing the Gromov product in $X$}
\label{figure:injecttriang}
\end{center}
\end{figure}
In this figure, all the edges of the inner triangle have length $1$.
Lemma
\ref{lemma:Tbottle} ensures that all the
edges leading from the inner to the outer triangle have length at most $2$.  A
computation then shows that $(x_i.x_j')\leq d((1,\nu),(1,\omega))+11$, and thus
$i([\phi])\neq i([\phi'])$.  
\end{proof}

\begin{definition}
A tree or quasi-tree $X$ is said to be \emph{bushy}, with
\emph{bushiness constant} $B$, if removing any metric ball of radius
$B$ from $X$ leaves a space with at least three unbounded path
components.  
\end{definition}
\begin{remark}
If two quasi-trees are quasi-isometric and one is bushy, then clearly
the other must also be bushy, though possibly with a different
bushiness constant.
\end{remark}
\begin{theorem}\label{th:pseudoquasibushy}
If $f\co G\to \R$ is a pseudocharacter which is not uniform, then $G$ admits a
cobounded quasi-action on a bushy tree.
\end{theorem}
\begin{proof}
By Corollary
\ref{cor:qi}, the space $X$ defined in \ref{definition:X} is
quasi-isometric to a tree $\Gamma$.  If $f$ is uniform, $E(f)$ contains at
least three points.
By Lemma \ref{lemma:injective} this implies that $\partial \Gamma
\cong \partial X$ 
contains
at least three points.  Let $\gamma_1$, $\gamma_2$, and $\gamma_3$ be
geodesic rays in $\Gamma$ starting at some fixed point $p\in\Gamma$
and tending to these three points in $\partial \Gamma$.  As these
points are distinct, there is some $R>0$ so that
$\gamma_1[R,\infty)\cap\gamma_2[R,\infty)\cap\gamma_3[R,\infty)$ is
empty.  Thus removing an $R$--ball centered at $p$ from $\Gamma$ leaves
a space with at least three unbounded path components.  Since $G$
quasi-acts coboundedly on $\Gamma$ (Corollary \ref{cor:qi}), there is
some constant $B$ so that removing \emph{any} $B$--ball leaves a space
with at least three unbounded path components.  
\end{proof}
\begin{remark}
Recall that the pseudocharacter $f$ gives rise to a $(1,\cobf)$--quasi-action on
$\R$.
The function $\varpi\co X\to \R$ given by $\varpi((g,\tau))=f(g)+f(\tau)$
is coarsely equivariant in an obvious sense.  Say that a $G$--quasi-tree 
(or tree with a $G$--quasi-action)
$\Lambda$ is \emph{maximal with respect to $f$} if $\Lambda$ admits a coarsely
equivariant map to $\R$ and if 
whenever $q\co\Lambda'\to\Lambda$ is a coarsely
equivariant coarsely surjective map from another $G$--quasi-tree or tree with a
$G$--quasi-action $\Lambda'$, then $q$ must be a quasi-conjugacy.  
It might be interesting to investigate the following questions:  
Do maximal quasi-trees exist?  Under what conditions is our quasi-tree $X$
maximal?  
\end{remark}

\subsection{Space of Pseudocharacters}
In this subsection we show that if there is a bushy pseudocharacter on $G$, then
the space of pseudocharacters on $G$ is actually infinite-dimensional.  We
first recall the terminology and main result of \cite{bestvinafujiwara:mcg}.
We consider the action of $G$ on the quasi-tree $X$ defined in the last
subsection.  Since $X$ is quasi-isometric to a $0$--hyperbolic space (a tree)
it is $\delta$--hyperbolic for some $\delta$.  The following definitions make
sense whenever $G$ is a group acting on a $\delta$--hyperbolic graph $X$ (see
\cite{bestvinafujiwara:mcg} for more details --  Definitions
\ref{def:hyperbolic}--\ref{def:independent}
are quoted
nearly verbatim from there).

\begin{definition}\label{def:hyperbolic}
Call an isometry $g$  of $X$ \emph{hyperbolic} if it admits a 
$(K,L)$--quasi-axis for some 
$K$, and $L$.  That is, there is a bi-infinite $(K,L)$--quasi-geodesic which
is mapped to itself by a nontrivial translation.  
This quasi-geodesic is said to be given the
\emph{$g$--orientation} if it is oriented so that $g$ acts as a positive
translation.  Note that any two $(K,L)$--quasi-axes of $g$ are within some
universal $B=B(\delta,K,L)$ of one another (by an elementary extension of Lemma 
\ref{lemma:quasigeodesic}), and any sufficiently long 
$(K,L)$--quasi-geodesic 
arc in a $B$--neighborhood of a quasi-axis for $g$ inherits a
natural $g$--orientation.
\end{definition}
\begin{definition}
If $g_1$ and $g_2$ are hyperbolic elements of $G$, write $g_1\sim g_2$ if for an
arbitrarily long segment $J$ in a $(K,L)$--quasi-axis for $g_1$ there is a
$g\in G$ such that $g(J)$ is within $B(\delta,K,L)$ of a $(K,L)$--quasi-axis of
$g_2$ and $g\co J\to g(J)$ is orientation-preserving with respect to the
$g_2$--orientation on $g(J)$.
\end{definition}

\begin{definition}\label{def:independent}
Two hyperbolic isometries $g_1$ and $g_2$ are said to be \emph{independent} if
their quasi-axes do not contain rays which are a finite Hausdorff distance
apart.  Equivalently %??%
the fixed point sets of $g_1$ and $g_2$ in $\partial X$ are disjoint.  An action
is \emph{nonelementary} if there are group elements which act as independent
hyperbolic isometries.
\end{definition}

\begin{definition}
A \emph{Bestvina-Fujiwara} action is a nonelementary action of a group $G$ on
 a hyperbolic graph $X$ so that there exist independent $g_1$, $g_2\in G$ so
that $g_1\not\sim g_2$.
\end{definition}
\begin{theorem}\label{thm:bestvinafujiwara}
{\rm\cite{bestvinafujiwara:mcg}}\qua  If $G$ admits a Bestvina-Fujiwara action, then 
$H_b^2(G;\R)$ and the space of pseudocharacters on $G$ both have dimension 
equal to $|\R|$.  
\end{theorem}

\begin{proposition}\label{prop:bestvinafujiwara}
If $f\co G\to \R$ is a bushy pseudocharacter, then the action on $X$ described
in Definition \ref{definition:X} is a Bestvina-Fujiwara action.
\end{proposition}
\begin{proof}
First note that if $g\in G$ and $f(g) > 0$, then $g$ acts as a hyperbolic
isometry of $X$.  Indeed, let $x_0=(g_0,\tau_0)$ be a vertex of $X$, and let
$x_n=g^n x_0 = (g^n g_0,\tau_0)$.  We choose a constant speed 
geodesic path from $x_0$ to $x_1$
and translate it to get a map $\gamma\co\R\to X$ so that $\gamma(n) = x_n$ for
all $n\in \Z$.  For $s$, $t\in \R$ we clearly have
$d(\gamma(s),\gamma(t))\leq D|s-t|$, where $D$ is the distance between $x_0$ and
$x_1$.  

Let $\overline{f}(g,\tau)=f(g)+f(\tau)$.  Suppose $x$, $x'$ are vertices of $X$ which
are connected by an edge.  It is straightforward to see that
$|\overline{f}(x)-\overline{f}(x')|\leq 2+4\cobf$.  This gives an upper bound for the 
gradient of $\overline{f}$ on $X$.  Since $|\overline{f}(x_m)-\overline{f}(x_n)|\geq
|m-n|f(g)-\cobf$, we get $d(\gamma(n),\gamma(m))\geq (|m-n|f(g)-\cobf)/( 2+4\cobf)$.
For arbitrary $s$, $t\in \R$ we have
\[d(\gamma(s),\gamma(t))>\frac{f(g)}{2+4\cobf}|s-t|-\left(\frac{\cobf}{2+4\cobf}+2K\right).\]
Thus $\gamma$ is a quasi-axis for $g$. 

\begin{claim}
If $g_1\sim g_2$, and $\sigma(g_1)\neq 0$, then 
$\sigma(g_2)=\sigma(g_1)$, where $\sigma(g)$ is the sign of $f(g)$ as in
Definition \ref{def:sigma}.
\end{claim}
\begin{proof}
Assume for simplicity that $\sigma(g_1)=1$.  It is sufficient to show that
$f(g_2^N)>0$ for some $N\geq 1$.  For $i\in\{1,2\}$, let $\gamma_i$ be a
$(K,L)$--quasi-axis for $g_i$, parameterized so that for some point $x_i$
$\gamma(n)=g_i^n(x_i)$.  

For $N>0$ there is an $h=h_N$ in $G$ so that $h(\gamma_2[0,N])$ is in a
$B(\delta,K,L)$--neighborhood of $\gamma_1$, where $\delta$ is the thinness
constant for $X$.  Furthermore, if $N$ is large enough, and $\gamma_1(p)$ and
$\gamma_1(q)$ are the closest elements of $\gamma_1(\Z)$ to $h\gamma_2(0)$ and
$h\gamma_2(N)$ respectively, then $q>p$.  We assume that $N$ is at least this
large.  Indeed, by choosing $N$ large enough, we can ensure that $q-p$ is as
large as we like.  

Since $|\overline{f}(g_1^n(x_1))-(nf(g_1)+\overline{f}(x_1))|$ is
bounded, we can therefore ensure that
$\overline{f}(\gamma(q))-\overline{f}(\gamma(p))$ is arbitrarily
large.  As the endpoints of $h\gamma_2[0,N]$ are at most
$B(\delta,K,L)+d(x_1,g_1x_1)$ from the points $\gamma_1(p)$ and
$\gamma_2(q)$, and the gradient of $\overline{f}$ is bounded, it
follows that we can choose $N$ to make
$\overline{f}(h\gamma_2(N))-\overline{f}(h\gamma_2(0))$ very large.
Thus we can make
$\overline{f}(\gamma_2(N))-\overline{f}(\gamma_1(0))$, and finally
$f(g_2^N)$, as large as we like.  In particular, we may find $N$ so
that $f(g_2^N)$ is positive.
\end{proof}

In Theorem \ref{th:freegroup}, we showed that there are elements $g$ and $g'$,
with $\sigma(g)=\sigma(g')=1$, so that $g$ and $g'$ have disjoint fixed point
sets in $E(f)$.  Thus $g^{-1}$ and $g'$ act independently and hyperbolically 
on $X$.
Furthermore, $\sigma(g')\neq \sigma(g^{-1})$, and so by the claim,
$g'\not\sim g^{-1}$.
\end{proof}

The following theorem is an immediate corollary of Proposition
\ref{prop:bestvinafujiwara} and Theorem \ref{thm:bestvinafujiwara}:
\begin{theorem}\label{th:infdim}
If $G$ admits a single bushy pseudocharacter, then $H_b^2(G;\R)$ and the 
space of pseudocharacters on $G$ both have dimension equal to $|\R|$.  
\end{theorem}

\section{Examples}\label{section:examples}

In \cite{msw:quasiactI}, it is shown that any quasi-action on a bounded valence
bushy tree is quasi-conjugate to an action on a (possibly different) bounded
valence bushy tree.  One might conjecture that
some kind of analogous statement holds
in the unbounded valence case.  An immediate obstacle to such a conjecture is
that any isometric
action on an $\R$--tree gives rise to a quasi-action on a simplicial
tree.  One way to see this is that $\R$--trees clearly
satisfy the bottleneck property of Theorem \ref{th:bottleneck},
and thus are quasi-isometric to simplicial trees.  One might still ask the
following question:
\begin{question}\label{question:naive}
Is every quasi-action on a bushy tree by a finitely presented group 
quasi-conjugate to an action on an $\R$--tree?
\end{question}

Kevin Whyte pointed out the following 
simple example after seeing an earlier version of
this paper: 
\begin{example}
Let $F$ be the graph with vertex set equal to 
$\overline{\Q}=\Q\cup\{\infty\}$ and so that two vertices $\frac{p}{q}$ and
$\frac{r}{s}$ are connected with an edge whenever $ps-qr=\pm 1$ (Formally we
think of $\infty$ as $\frac{1}{0}$.).  This is usually called the Farey graph.
 The group $PSL(2,\Z)$ acts on the vertices by M\"{o}bius transformations, and
preserves the edge condition.  It is not hard to see that $F$ satisfies the
bottleneck condition of Theorem \ref{th:bottleneck}, and is thus a quasi-tree.
In fact, $F$ is quasi-isometric to an infinite valence tree.  The action of
$PSL(2,\Z)$ on $F$ thus induces a quasi-action on an infinite valence tree, 
via Proposition \ref{prop:equiv}.  
\end{example}
\begin{proposition}
The action of $PSL(2,\Z)$ on the Farey graph is not quasi-conjugate to an
isometric action on any $\R$--tree.
\end{proposition}
\begin{proof}[Sketch proof]
Suppose $PSL(2,\Z)$ acts isometrically on an $\R$--tree $\Lambda$, and that this action is quasi-conjugate to the action on the Farey graph $F$.  
Note that $PSL(2,\Z)$ is generated by the elements
$A=\left(\begin{array}{cc} 1 & 1 \\ 0 & 1
                                                \end{array}\right)/\{\pm I\}$ 
and 
$B=\left(\begin{array}{cc} 1 & 0 \\ -1 & 1
                                                \end{array}\right)/\{\pm I\}$.
Both $A$ and $B$ must fix points in $\Lambda$, as they fix points in $F$. ($A$
fixes $\frac{1}{0}$ and $B$ fixes $\frac{0}{1}$.)
The element $(AB)^{-1}=\left(\begin{array}{cc} 1 & -1 \\ 1 & 0
                                                \end{array}\right)/\{\pm I\}$ 
has
finite order, and so it must also 
it must also fix a point
in $\Lambda$.  It then follows (eg from \cite[Corollary~1 on p64]{serre:trees}) that $A$ and $B$ must have a common fixed point.
  Since $A$ and $B$ generate, this
implies that $PSL(2,\Z)$ fixes a point in $\Lambda$.  Thus every orbit in 
$\Lambda$ has
finite diameter.  This would imply that every orbit in $F$ has finite diameter,
which is easily seen to be false.    
\end{proof}

Thus the answer to Question \ref{question:naive} is no.  Of course
$PSL(2,\Z)$ admits a cocompact action on the (infinite diameter)
 Bass-Serre tree coming
from the splitting $PSL(2,\Z)\cong \Z_2\ast\Z_3$. 
 One might still ask the
following question:

\begin{question}
Does every group with a cobounded quasi-action on a bushy tree also act
nontrivially and isometrically on some tree?
\end{question}

We will show that the answer to this question is also no.
 
\begin{figure}[ht!]
\begin{center}
\input{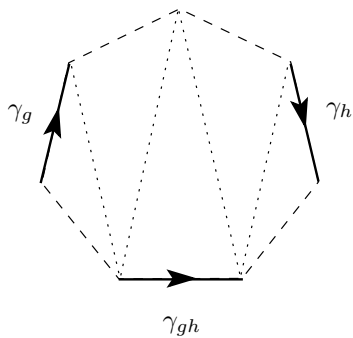}
\caption{Pants cut open and triangulated with straight triangles}
\label{figure:heptagon}
\end{center}
\end{figure}

Suppose that $M$ is a closed Riemannian 
manifold with all sectional curvatures $\leq -1$.  Here
is one way to generate examples of pseudocharacters on $\pi_1(M)$.  Let $\omega$
be any $1$--form on $M$.  If $g\in \pi_1(M)$ we let $\gamma_g$ be the unique
closed geodesic in its free homotopy class.  Let $f_\omega\co\pi_1(M)\to\R$
 be given by 
\begin{equation}\label{eq:formdef}
f_\omega(g)=\int_{\gamma_g}\omega.
\end{equation}
We claim that $f_\omega$ is a pseudocharacter on $G$.  It is clear from the
definition that $f_\omega$ is conjugacy invariant and a homomorphism on each
cyclic subgroup.  To see it is a coarse homomorphism on $G$, we compute for
$g$, $h\in\pi_1(M)$,
\[\delta f_\omega(g,h)=f_\omega(gh)-f_\omega(g)-f_\omega(h)\]
\[=\int_{\gamma_{gh}\cup -\gamma_g\cup -\gamma_h}\omega = \int_F d\omega\]
where $F$ is a (not necessarily embedded) pair of pants in $M$ like the one 
shown in Figure
\ref{figure:pants}.
\begin{figure}[ht!]
\begin{center}
\input{pants.pstex_t}
\caption{Pants}
\label{figure:pants}
\end{center}
\end{figure}
The quantity $\int_F d\omega$ only depends on $F$ up to homotopy in $M$ so we
may assume that $F$ is triangulated as in Figure \ref{figure:heptagon} and that
each triangle of $F$ has been straightened.  This is to say that each edge of
the triangulation has been made geodesic, and each face is a union of geodesics
issuing from one vertex and terminating at the opposite edge.
One may show using the Gauss-Bonnet theorem that the area of any straight 
triangle in $M$ is $\leq \pi$.
 Thus the area of $F$ is $\leq 5\pi$ after straightening.  
Since $M$ is compact, $d\omega$ is bounded, and so $\int_F d\omega$
is bounded above by some multiple of the area of $F$.  Thus
$\|\delta f_\omega\|<\infty$.  

We would like to use the preceding construction to obtain a bushy
pseudocharacter on a group with no nontrivial action on any tree.
Thurston showed in \cite{thurston:79} how to obtain (via Dehn filling)
many negatively
curved three-manifolds whose fundamental groups cannot act
nontrivially on any tree.
\begin{proposition}
All but finitely many fillings of the figure eight knot complement have 
fundamental groups which
admit bushy pseudocharacters.  
\end{proposition}
\begin{proof}
Let $M$ be the complement of the figure eight knot in $S^3$.
In \cite{thurston:79}, Thurston showed that $M$
admits a complete hyperbolic metric of finite volume.
In \cite{bart:surfacegroups}, Bart shows that $M$
contains
a closed, immersed, totally geodesic surface $\Sigma$ which remains 
$\pi_1$--injective
after all but at most thirteen fillings.  
The main tool used in Bart's proof is
the Gromov-Thurston $2\pi$ Theorem, 
in which an explicit negatively curved
metric is constructed on the filled manifold \cite{bleilerhodgson:twopi}.

Let $M(\gamma)$ be one of the fillings which can be given a negatively curved
metric, and let $G=\pi_1 M(\gamma)$.  
We assume that $M(\gamma)$ is endowed with the Riemannian metric given by the 
$2\pi$ theorem.
Outside of a neighborhood of the core curve of the filling solid torus,
this metric is isometric to the hyperbolic metric on $M$ with a neighborhood of
the cusp removed.  
Thus the surface $\Sigma$ remains totally geodesic in $M(\gamma)$, and
does not intersect the filling solid torus.
Furthermore, the core curve $c$ of the filling solid torus is a closed
geodesic.  

Let $\omega$ be a one-form on $M(\gamma)$, supported inside the filling
solid torus, so that $\int_c\omega > 0$, and define
$f_\omega\co\pi_1(M(\gamma))\to\R$ as in equation \ref{eq:formdef}.  We will show
that $f_\omega$ is a bushy pseudocharacter. 

We fix a generating set $S$ for $G$ and consider the Cayley graph
$\Gamma=\Gamma(G,S)$.  Since $M(\gamma)$ is compact, $\Gamma$ is
quasi-isometric to the universal cover $\widetilde{M}(\gamma)$.  Since
$M(\gamma)$ has a negatively curved Riemannian metric, $\Gamma$ must
be negatively curved in the sense of Gromov;  further, the Gromov
boundary $\partial \Gamma$ may be identified with the $2$--sphere which
is the visual boundary of
$\widetilde{M}(\gamma)$ \cite{bridhaef:book}.  

The inclusion of $\Sigma$ into $M(\gamma)$
induces an inclusion of the surface group $F=\pi_1(\Sigma)$ into
$G$; we note that some neighborhood $N(F)$ of $F$ separates $\Gamma$
into two unbounded complementary components.  As $\Sigma$ is totally
geodesic in $M(\gamma)$, the Gromov boundary $\partial F\cong S^1$
embeds in $\partial G$ and cuts $\partial G$ into two open disks $D_1$
and $D_2$.  

If $g$ is any element conjugate into $F$, then its geodesic
representative actually lies in $\Sigma$, and so $f_\omega(g)=0$.

Let $g_c\in G$ be some
element whose geodesic representative is $c$.  Then $[g_c^\infty]\in
E(f)^+$, and $\{g_c^n\}_{n\in\N}$ limits to some point $e_c\in\partial
G$.  The orbit of this point under the action of $G$ on $\partial G$
is dense \cite{gromov:wordhyperbolic}, so we may find $h_1$
and $h_2$ in
$G$ so that $h_ie_c\in D_i$ for $i\in\{1,2\}$.

For each $i\in\{1,2\}$, let
$\phi_i\co\R_+\to\Gamma$ be a constant speed path with $\phi_i(0)=h_i$
and $\phi_i(n)=h_i g_c^n$.  Then $\lim_{t\to\infty} \phi_i(t) = h_i
e_c$, and $[\phi_i]\in E(f)^+$ for each $i$.  

We claim that $[\phi_1]\neq[\phi_2]$ in $E(f)$.  Indeed, if
$[\phi_1]$ and $[\phi_2]$ are equal, then (applying Definition
\ref{def:v1}) there are connecting
paths $\delta\co[0,1]\to\Gamma$ between $\phi_1$ and $\phi_2$ so that
$\sup(f\circ\delta([0,1]))$ is large, but the diameter of
$f\circ\delta([0,1])$ is small.  Because the paths $\phi_1$
and $\phi_2$ go into separate components of the complement of $N(F)$,
these connecting paths must pass through $N(F)$ (Figure \ref{figure:fomegaisbushy}).
\begin{figure}[ht!]
\begin{center}
\input{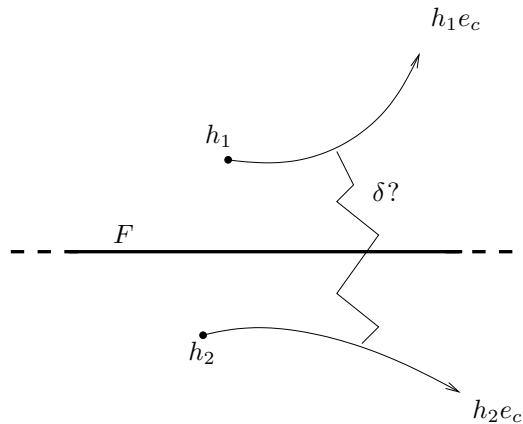}
\caption{$\phi_1$ and $\phi_2$ cannot represent the same element of $E(f)$.}
\label{figure:fomegaisbushy}
\end{center}
\end{figure}
Since $f_\omega$ is
zero on $F$, it is bounded on $N(F)$, and these connecting paths
cannot exist.  

A similar argument shows that we can find distinct elements of
$E(f_\omega)^-$.
\end{proof}

\begin{remark}
The proof of the preceding proposition can be applied to fillings of any
hyperbolic $3$--manifold with an immersed closed incompressible surface, so long
as the surface is quasi-Fuchsian.  
\end{remark}

It is proved in \cite{morganshalen:degenerationsIII} 
that if a $3$--manifold group acts on
an $\R$--tree non-trivially, then it must split as either an amalgamated free
product or as an HNN extension.  By a standard argument, this implies that the
 $3$--manifold is either reducible or Haken.  Since all but finitely many
fillings of the figure eight knot complement are non-Haken,
 we have the following
corollary:
\begin{corollary}\label{cor:counter}
There are infinitely many closed $3$--manifold groups which quasi-act 
coboundedly on bushy trees but which admit no
nontrivial isometric action on any $\R$--tree.
\end{corollary}

\begin{remark}
In a future article \cite{manning:qfa} we will look more closely at
rigidity and irrigidity of 
group quasi-actions on trees, and in particular
give examples of groups (for example $SL(n,\Z)$ for $n>2$) which
do not quasi-act coboundedly on any (infinite diameter) tree.  
\end{remark}

\end{document}